\documentclass{article}%
\usepackage{amsmath}
\usepackage{graphicx}
\usepackage{amsfonts}
\usepackage{amssymb}%
\providecommand{\U}[1]{\protect\rule{.1in}{.1in}}
\providecommand{\U}[1]{\protect\rule{.1in}{.1in}}
\newtheorem{theorem}{Theorem}
\newtheorem{acknowledgement}[theorem]{Acknowledgement}

\newtheorem{definition}[theorem]{Definition}
\newtheorem{example}[theorem]{Example}

\newtheorem{remark}[theorem]{Remark}

\begin{document}

\title{Towards a unified theory of Sobolev inequalities}
\author{Joaquim Mart\'{\i}n\thanks{Partially supported in part by Grants
MTM2010-14946, MTM-2010-16232}\\Department of Mathematics\\Universitat Aut\`{o}noma de Barcelona\\jmartin@mat.uab.cat
\and Mario Milman\thanks{ Partially supported by a grant from the Simons Foundation
(\#207929 to Mario Milman).}\\Department of Mathematics\\Florida Atlantic University\\mario.milman@gmail.com\\https://sites.google.com/site/mariomilman/}
\date{}
\maketitle

\tableofcontents

\begin{abstract}
We discuss our work on pointwise inequalities for the gradient which are
connected with the isoperimetric profile associated to a given geometry. We
show how they can be used to unify certain aspects of the theory of Sobolev
inequalities. In particular, we discuss our recent papers on fractional order
inequalities, Coulhon type inequalities, transference and dimensionless
inequalities and our forthcoming work on sharp higher order Sobolev
inequalities that can be obtained by iteration.

\end{abstract}

\section{Introduction}

In this expository note we survey some of our previous work towards a unified
theory of Sobolev inequalities using pointwise rearrangement inequalities. The
presentation includes our recent results \cite{mamicoulhon} on Coulhon
inequalities (cf. \cite{cou1}, \cite{bakr}), integral transference and
dimensionless inequalities \cite{mamidim}, fractional inequalities
\cite{mamiaster}, extrapolation and self improvement \cite{mamicon}. We also
include a brief section on how iterations of our estimates lead to sharp
higher order pointwise inequalities \cite{mamiitera}. We hope that this
attempt to summarize and organize some of the material, together with the
inclusion of our motivation and details on the methods used, could be of some interest.

In conclusion we should mention three important lessons we learned from Gian
Carlo Rota \cite{rota}:

\begin{center}
\textbf{1}\textit{. Publish the same results several times \textbf{2.} Do not
worry about your mistakes and \textbf{3.} Write informative introductions.}
\end{center}

\section{Basic definitions and notation}

Before start our presentation let us recall some basic definitions and set up
the notation.

Let $\left(  \Omega,\mu\right)  $ be a measure space $\left(  \Omega
,\mu\right)  ;$ for a measurable function $u:\Omega\rightarrow\mathbb{R},$ the
\textbf{distribution function} of $u$ is given by
\[
\mu_{u}(t)=\mu\{x\in{\Omega}:\left\vert u(x)\right\vert >t\}\text{
\ \ \ \ }(t\geq0).
\]
The \textbf{decreasing rearrangement} of a function $u$ is the
right-continuous non-increasing function from $[0,\mu(\Omega))$ into
$\mathbb{R}$ which is equimeasurable with $u.$ It can be defined by the
formula%
\[
u_{\mu}^{\ast}(s)=\inf\{t\geq0:\mu_{u}(t)\leq s\},\text{ \ }s\in\lbrack
0,\mu(\Omega)).
\]
The maximal average $u_{\mu}^{\ast\ast}(t)$ is defined by
\[
u_{\mu}^{\ast\ast}(t)=\frac{1}{t}\int_{0}^{t}u_{\mu}^{\ast}(s)ds.
\]
When the measure is clear from the context, or when we are dealing with
Lebesgue measure, we may simply write $u^{\ast}$ and $u^{\ast\ast}$, etc.

Let $X=X({\Omega})$ be a Banach function space on $({\Omega},\mu)$, with the
Fatou property\footnote{This means that if $f_{n}\geq0,$ and $f_{n}\uparrow
f,$ then $\left\Vert f_{n}\right\Vert _{X}\uparrow\left\Vert f\right\Vert
_{X}$ (i.e. Fatou's Lemma holds in the $X$ norm).}. We shall say that $X$ is a
\textbf{rearrangement-invariant} (r.i.) space, if $g\in X$ implies that all
$\mu-$measurable functions $f$ with $f_{\mu}^{\ast}=g_{\mu}^{\ast},$ also
belong to $X$ and moreover, $\Vert f\Vert_{X}=\Vert g\Vert_{X}$. The
functional $\Vert\cdot\Vert_{X}$ \ will be called a rearrangement invariant
norm. Typical examples of r.i. spaces are the $L^{p}$-spaces, Orlicz spaces,
Lorentz spaces, Marcinkiewicz spaces, etc.

Let $X({\Omega})$ be a r.i. space, then there exists a r.i. space $\bar
{X}=\bar{X}(0,\mu(\Omega))$ on $\left(  (0,\mu(\Omega)),m\right)  $, ($m$
denotes the Lebesgue measure on the interval $(0,\mu(\Omega))$) such that%
\begin{equation}
\Vert f\Vert_{X({\Omega})}=\Vert f_{\mu}^{\ast}\Vert_{\bar{X}(0,\mu(\Omega))}.
\label{repppp}%
\end{equation}
$\bar{X}$ is called the \textbf{representation space} of $X({\Omega})$ which
is unique if $\mu(\Omega)$ is finite.) We refer the reader to \cite{BS} for a
complete treatment of the theory of (r.i.) spaces.)

As usual, the symbol $f\simeq g$ will indicate the existence of a universal
constant $c>0$ (independent of all parameters involved) so that $(1/c)f\leq
g\leq c\,f$, while the symbol $f\preceq g$ means that $f\leq c\,g,$ finally
given $1\leq p\leq\infty$, $p^{\prime}=p/(p-1)$.

\section{Two Examples}

We state two simple examples that illustrate the difficulties in deriving a
unified theory of Sobolev inequalities.

\begin{example}
In dimension 1 the Sobolev embedding theorem is connected with the fundamental
theorem of calculus: If $f\ \ $is a Lip function with compact support, then%
\[
f(y)=\int_{-\infty}^{y}f^{\prime}(s)ds
\]
and therefore\footnote{actually using the available cancellation, we have
$2f(y)=\int_{-\infty}^{y}f^{\prime}(s)ds+\int_{y}^{\infty}f^{\prime}(s)ds$,
therefore the constant of the embedding can be improved,%
\[
\left\Vert f\right\Vert _{\infty}\leq\frac{1}{2}\left\Vert f^{\prime
}\right\Vert _{1}.
\]
}%
\[
\left\Vert f\right\Vert _{\infty}\leq\left\Vert f^{\prime}\right\Vert _{1}.
\]
In dimension $n$ the corresponding estimate does not hold. In particular,
there are functions $f\in W^{1,n}(\mathbb{R}^{n})$ that are not bounded.
Instead, we have somewhat weaker results: (Sobolev)%
\[
\left\Vert f\right\Vert _{p_{n}}\leq c_{n,p}\left\Vert \nabla f\right\Vert
_{p},f\in C_{0}^{\infty}(\mathbb{R}^{n}),\text{ with }\frac{1}{p_{n}}=\frac
{1}{p}-\frac{1}{n},\text{ }1<p<n,
\]
with $c_{n}\rightarrow\infty,$ when $p\rightarrow n.$ The inequality is also
true when $p=1$ (Gagliardo-Nirenberg): there exists $c_{n}$ such that%
\[
\left\Vert f\right\Vert _{n^{\prime}}\leq c_{n}\left\Vert \nabla f\right\Vert
_{1},\text{ }f\in C_{0}^{\infty}(\mathbb{R}^{n}).
\]
One way to complete the picture for $p=n,$ is by extrapolation: Trudinger
\cite{tr} showed by extrapolation that if $\Omega$ is an open set in
$\mathbb{R}^{n}$ with $\left\vert \Omega\right\vert <\infty$, then
\[
W_{0}^{1,n}(\Omega)\subset e^{L^{n^{\prime}}}.
\]

\end{example}

This example, in particular, illustrates the crucial role of dimension in the
Euclidean Sobolev inequalities.

As a counterpart, consider the following

\begin{example}
Consider $\mathbb{R}$ but now equipped with one dimensional Gaussian measure
$\gamma_{1}$. In this case the gain of integrability is logarithmic
(\textquotedblleft logarithmic type Sobolev inequalities\textquotedblright)
\[
\left\Vert f\right\Vert _{L(LogL)^{1/2}(\gamma_{1})}\leq c(\left\Vert
\left\vert \nabla f\right\vert \right\Vert _{L^{1}(\gamma_{1})}+\left\Vert
f\right\Vert _{L^{1}(\gamma_{1})}).
\]
Similarly, for $p>1$ we have%
\[
\left\Vert f\right\Vert _{L^{p}(LogL)^{p/2}(\gamma_{1})}\leq c(\left\Vert
\left\vert \nabla f\right\vert \right\Vert _{L^{p}(\gamma_{1})}+\left\Vert
f\right\Vert _{L^{p}(\gamma_{1})}).
\]
Moreover, the gain of integrability does not change with the dimension! For
example, for $\mathbb{R}^{n}$ equipped with Gaussian measure $\gamma_{n}$, we
have%
\[
\left\Vert f\right\Vert _{L(LogL)^{1/2}(\gamma_{n})}\leq c(\left\Vert
\left\vert \nabla f\right\vert \right\Vert _{L^{1}(\gamma_{n})}+\left\Vert
f\right\Vert _{L^{1}(\gamma_{n})}).
\]
In other words, while in the Euclidean geometry the spaces involved in the
inequalities depend on the dimension, and the gain of integrability is
naturally measured with powers, i.e. using $L^{p}$ or $L(p,q)$ spaces, in the
Gaussian world the gain of integrability is logarithmic, independent of the
dimension, and the right spaces one needs to measure the gain of integrability
are logarithmic Orlicz type spaces. The celebrated logarithmic Sobolev
inequalities of Gross (cf. \cite{gro}) are among the most prominent examples
of dimensionless inequalities in the Gaussian world (cf. \cite{gro}).
\end{example}

Traditionally this state of affairs has led to different methods/theories to
deal with Euclidean or Gaussian Sobolev inequalities, or more generally,
Sobolev inequalities in other geometries, e.g. Euclidean domains with measures
of the form $w(x)dx,$ Riemannian manifolds, or more generally, metric measure
spaces. In particular, what is the role of dimension? What are the function
spaces one needs to use to measure the integrability gains?

While some of the differences are unavoidable, one wonders if it is possible
to unify at least some aspects of these disparate theories and thus, maybe,
provide a better understanding. Recently there has been progress in this
direction from several different directions by a number of authors. The
results presented in this note describe mainly of our work towards a unified
theory of Sobolev inequalities using pointwise inequalities on rearrangements;
as consequence the list of references is disproportionately tilted towards our
own work. Moreover, there is a huge literature on Sobolev inequalities, and
while the bibliography we have compiled is rather large, we must warn the
reader that we did not attempt to cover all the references, not even the
important references. We apologize in advance if your favorite paper/author is
not included in the bibliography, hopefully it should not be more than one
(reference) iteration away. In particular, we should explicitly mention our
debt to the pioneering\footnote{For a detailed presentation of Maz'ya's
remarkable early work we refer to \cite{hedb}.} work of V. Maz'ya which has
greatly influenced our view point of the subject (cf. Section \ref{sec:maz}
below). For more information, background and more comprehensive bibliographies
we refer to \cite{adams}, \cite{edmund}, \cite{HK2}, \cite{leoni},
\cite{mamiadv}, \cite{maz}, \cite{pick}, \cite{ra}, \cite{saloff coste}, and
\cite{tal}.

\section{Rearrangement invariant spaces and the Brezis-Wainger-Hansson
embedding}

One early difficulty in trying to develop any type of unified theory is that,
even in the Euclidean case, the limiting borderline case ($p=n)$ of the
Sobolev embedding apparently requires the use of a different scale of spaces,
e.g. an \textquotedblleft extrapolation space"...So we started our work trying
to understand the issues connected with the limiting inequalities. In a
convoluted way our efforts to understand the limiting cases eventually led us
to a better understanding on how to approach all the cases
simultaneously\footnote{This is a somewhat dissapointing turn of events for
the developers of general abstract theories studying limiting inequalities
(e.g. \cite{jm}) but our current understanding of Sobolev inequalities shows
that: (a) Sobolev inequalities self improve (cf. \cite{bakr}, \cite{haj},
\cite{mmp}, \cite{mamicon}) and (b) the extrapolations of Sobolev inequalities
take the form: \textquotedblleft one inequality" implies a family of
inequalities and in some cases \textquotedblleft one inequality implies all"
!(cf. \cite{bakr}, \cite{cwijami}, \cite{mamicon})). We also refer to the
forthcoming \cite{mamirubio} for a connection with extrapolation of Sobolev
inequalities \`{a} la Rubio de Francia.}! So it seems that this is a good
point where to start with the story.

An early result in this direction was obtained by Trudinger \cite{tr}, who
showed that for a domain $\Omega$ in $\mathbb{R}^{n},$ with $\left\vert
\Omega\right\vert <\infty,$%
\begin{equation}
W_{0}^{1,n}(\Omega)\subset e^{L^{n^{\prime}}}(\Omega). \label{oldcase}%
\end{equation}
Brezis-Wainger \cite{bw} improved this result using the rearrangement
inequalities of O'Neil \cite{on}, while Hansson \cite{hansson} obtained
similar results. Their result reads%
\begin{equation}
W_{0}^{1,n}(\Omega)\subset BWH(\Omega), \label{embed}%
\end{equation}
where if, say, $\left\vert \Omega\right\vert =1,$ then $BWH=$
(Brezis-Wainger-Hansson) is the space defined by%
\[
\left\Vert f\right\Vert _{BWH}=\left\{  \int_{0}^{1}\left(  f^{\ast\ast
}(s)\frac{1}{1+\log\frac{1}{s}}\right)  ^{n}\frac{ds}{s}\right\}
^{1/n}<\infty.
\]
Maz'ya (cf. \cite[(3.1.4) page 232]{maz} and the references therein) had
obtained earlier an inequality using his capacity theory that, in particular,
implies the embedding (\ref{embed}). The improvement over (\ref{oldcase}) is
given by the easily verified fact that%
\[
BWH(\Omega)\subset e^{L^{n^{\prime}}}(\Omega).
\]
On the other hand, O'Neil \cite{on} extended the original (one dimensional)
results of Hardy-Littlewood, as follows%
\[
W^{1,p}(\mathbb{R}^{n})\subset L(p_{n},p),\text{ }1<p<n,\text{ }\frac{1}%
{p_{n}}=\frac{1}{p}-\frac{1}{n},
\]
where the Lorentz $L(p,q)$ spaces, for $p<\infty,q\in\lbrack1,\infty]$ are
defined by%
\[
L(p,q)=\{f:\left\Vert f\right\Vert _{L(p,q)}=\left\{  \int_{0}^{\infty
}(f^{\ast\ast}(t)t^{1/p})^{q}\frac{dt}{t}\right\}  ^{1/q}<\infty\}.
\]
The results of O'Neil and Brezis-Wainger-Hansson, like the original one
dimensional results of Hardy-Littlewood, can be shown to be optimal within the
class of rearrangement invariant spaces: if $X$ is a rearrangement invariant
space then,%
\[
W^{1,p}\subset X\Rightarrow L(p_{n},p)\subset X,\text{ if }p<n,\text{ or
}BWH\subset X,\text{ if }p=n.
\]
As a consequence, the $L(p.q)$ spaces are not enough to describe the optimal
Sobolev inequalities (e.g. the case $p=n$ requires a different space). This
led to the introduction in the theory of Sobolev embeddings of the following
modification of the $L(p,q)$ spaces which, as we shall see, resolves this
difficulty (cf. \cite{bmr}, \cite{bmrind} and, as it turns out, in a different
way and less explicitly in \cite{tar} and \cite{mapi}). For a measure space we
define%
\[
L(\infty,\infty)=\{f:\left\Vert f\right\Vert _{L(\infty,\infty)}=\sup
_{t}(f^{\ast\ast}(t)-f^{\ast}(t))<\infty\}.
\]
This space was introduced by Bennett-DeVore-Sharpley \cite{bds} who in their
paper show that for functions defined on a cube, $L(\infty,\infty)$ is the
rearrangement invariant hull of $BMO,$ i.e. $L(\infty,\infty)$ is the smallest
possible space that contains all the rearrangements of functions in $BMO.$ One
should note here that the usual definition of $L(\infty,\infty)$ would give%
\[
\left\Vert f\right\Vert _{L(\infty,\infty)}=\sup_{t}f^{\ast\ast}(t)=\left\Vert
f\right\Vert _{L^{\infty}},
\]
while the space $L(\infty,\infty)$ that we have defined is bigger: we have
$BMO\subset L(\infty,\infty).$

More generally, note that if we formally attempt to define the $L(\infty,q)$
spaces using the classical definition the resulting spaces are trivial:%

\[
\int_{0}^{\infty}f^{\ast\ast}(s)^{q}\frac{ds}{s}<\infty\Rightarrow f=0.
\]
On the other hand, if we redefine the $L(\infty,q)$ spaces by means of
replacing $f^{\ast\ast}(t)$ by the oscillation $f^{\ast\ast}(t)-f^{\ast}(t),$
then the spaces defined by the condition%
\[
\left\Vert f\right\Vert _{L(\infty,q)}=\left\{  \int_{0}^{\left\vert
\Omega\right\vert }\left(  \left(  f^{\ast\ast}(t)-f^{\ast}(t)\right)
\right)  ^{q}\frac{dt}{t}\right\}  ^{1/q}<\infty,
\]
are not trivial, since the differentiation theorem provides us with a
cancellation at the origin. In \cite{bmr}, \cite{bmrind} the role of the
$L(\infty,q)$ spaces in the Sobolev embedding was observed. It was noted there
that these spaces were relevant in connection with an inequality implicit in
the paper by Alvino-Trombetti-Lions \cite{alv}: If $f$ is symmetrically
decreasing, then%
\begin{equation}
f^{\ast\ast}(t)-f^{\ast}(t)\leq c_{n}t^{1/n}\left\vert \nabla f\right\vert
^{\ast\ast}(t).\label{pll}%
\end{equation}
Now, if we let $f^{\circ}$ denote the symmetric rearrangement of $f$ $($cf.
\cite{kesavan}, \cite{leoni})$,$ then for smooth $f,$ the P\'{o}lya-Szeg\"{o}
principle can be formulated as (cf. \cite{mmp} and the references there to
earlier work by Fournier)
\[
\left\vert \nabla f^{\circ}\right\vert ^{\ast\ast}(t)\leq\left\vert \nabla
f\right\vert ^{\ast\ast}(t).
\]
Then, since $\left(  f^{\circ}\right)  ^{\ast}=f^{\ast},$ it follows that we
can eliminate the restriction for $f$ to be symmetrically decreasing,
therefore (\ref{pll}), indeed, holds for all smooth $f.$ Next, integrating
(\ref{pll}), we see that:%
\begin{equation}
\left\{  \int_{0}^{\left\vert \Omega\right\vert }\left(  \left(  f^{\ast\ast
}(t)-f^{\ast}(t)\right)  t^{1/p-1/n}\right)  ^{q}\frac{dt}{t}\right\}
^{1/q}\leq c_{n}\left\{  \int_{0}^{\left\vert \Omega\right\vert }\left(
t^{1/p}\left\vert \nabla f\right\vert ^{\ast\ast}(t)\right)  ^{q}\frac{dt}%
{t}\right\}  ^{1/q}.\label{limite}%
\end{equation}
The left hand side is equivalent to $\left\Vert f\right\Vert _{L(p_{n},q)}$
(cf. \cite{mami07}), and moreover, for $q=p$ we have
\[
\left\{  \int_{0}^{\left\vert \Omega\right\vert }\left(  t^{1/p}\left\vert
\nabla f\right\vert ^{\ast\ast}(t)\right)  ^{p}\frac{dt}{t}\right\}
^{1/p}\leq p^{\prime}\left\Vert \left\vert \nabla f\right\vert \right\Vert
_{p}.
\]
Thus, for $q=p<n$ we have recovered the classical Sobolev inequality.
Moreover, the inequality (\ref{limite}) is valid and makes sense in the
limiting case $p=n.$ In particular for $q=p=n,$ we have%
\[
\left\{  \int_{0}^{\left\vert \Omega\right\vert }\left(  \left(  f^{\ast\ast
}(t)-f^{\ast}(t)\right)  \right)  ^{n}\frac{dt}{t}\right\}  ^{1/n}\leq
c_{n}\left\Vert \nabla f\right\Vert _{n}.
\]
The condition that the left hand side of the previous inequality is finite
defines the space $L(\infty,n),$ and, moreover, we have (cf. \cite{bmrind})%
\[
L(\infty,n)\subset BWH.
\]
Since redefining the $L(p,q)$ spaces by means of replacing $f^{\ast\ast}$ by
$f^{\ast\ast}-f^{\ast}$ gives equivalent norms when the parameters are in the
usual range of the classical theory\footnote{For a detailed discussion on
equivalences between different Lorentz *norms* we refer to \cite{bmrind},
\cite{jm2}, \cite{mami07}.}, we now have a unified method to prove an
inequality that includes the problematic case $p=n,$%
\begin{equation}
\int_{0}^{\left\vert \Omega\right\vert }\left(  \left(  f^{\ast\ast
}(t)-f^{\ast}(t)\right)  t^{1/p-1/n}\right)  ^{q}\frac{dt}{t}\leq c_{n}%
\int_{0}^{\left\vert \Omega\right\vert }\left(  t^{1/p}\left\vert \nabla
f\right\vert ^{\ast\ast}(t)\right)  ^{q}\frac{dt}{t}.\label{flaca}%
\end{equation}
In fact, let us note that we can do this even if the measure is infinite,%
\[
\left\{  \int_{0}^{\infty}\left(  \left(  f^{\ast\ast}(t)-f^{\ast}(t)\right)
\right)  ^{n}\frac{dt}{t}\right\}  ^{1/n}\leq c_{n}\left\Vert \left\vert
\nabla f\right\vert \right\Vert _{n}.
\]
In this notation, the original Hardy-Littlewood-O'Neil program\footnote{The
proof to deal with the case $p=1$ is slightly different and hinges on a
variant of ( \ref{pll}), namely (cf. \cite{mmp})
\[
\int_{0}^{t}(f^{\ast\ast}(s)-f^{\ast}(s))s^{1/n}\frac{ds}{s}\leq c\int_{0}%
^{t}(\left\vert \nabla f\right\vert ^{\ast}(s))ds.
\]
} has been completed\footnote{Actually (\ref{flaca}) also makes sense, and
gives sharp results, when $p>n$ (cf. \cite[Chapter 9]{mamiaster}).}:%
\[
\left\Vert f\right\Vert _{L(p_{n},p)}\leq c\left\Vert \left\vert \nabla
f\right\vert \right\Vert _{L^{p}},\text{ }1\leq p\leq n.
\]
The improvement on Brezis-Wainger-Hansson is possible due to the fact that the
class $L(\infty,n)$ is not a linear space.

A posteriori, we also understood that an inequality obtained by Kolyada on the
unit cube \cite{kolyada}, that is exactly like (\ref{pll}), except that the
oscillation condition on the left hand side is given in terms of $f^{\ast
}(t)-f^{\ast}(2t)$ , could also have been used for the Euclidean inequalities
above (cf. \cite{perez}). Moreover, Tartar \cite{tar} (cf. also \cite{mapi})
earlier than \cite{bmrind}, but apparently after \cite{bmr}, had also
obtained, using a different approach based on truncations, a discrete version
of a result that, with some work, can be seen to be equivalent to
(\ref{flaca}). This was not clearly understood at the time we wrote
\cite{bmrind}.

From the point of view of the development of our program, we draw the
following from this section: Redefining the target spaces using oscillations,
and using the pointwise inequality (\ref{pll}), we could treat all the cases
of the classical Sobolev inequalities in an optimal unified manner.

To proceed further with our program of understanding more general geometries
the question we faced next was: What would be a substitute for (\ref{pll})
when dealing with other geometries? Since the oscillation condition
$f^{\ast\ast}(t)-f^{\ast}(t)$ is a general construct and, in fact, can be
understood from an approximation point of view (cf. \cite{jm2}), we
concentrated our efforts on understanding better the right hand side of the inequality.

\section{Martin-Milman-Pustylnik meet Maz'ya\label{sec:maz}}

The Gagliardo-Nirenberg inequality%
\begin{equation}
\left\Vert f\right\Vert _{L^{n^{\prime}}}\leq c\left\Vert \left\vert \nabla
f\right\vert \right\Vert _{L^{1}}, f\in C_{0}^{\infty}(\mathbb{R}^{n}),
\label{gagliardoN}%
\end{equation}
is well known to be equivalent to the isoperimetric inequality (cf. Maz'ya
\cite{maz} and the references therein.) Moreover, using the chain rule and the
scale properties of the $L^{p}$ spaces (here scale = H\"{o}lder's inequality)
one can readily see that (\ref{gagliardoN}) implies \textquotedblleft
all\textquotedblright\ the classical Sobolev inequalities: simply consider
positive $f$ and apply (\ref{gagliardoN}) to $f^{\alpha}$ for suitable
$\alpha$ combined with the chain rule and H\"{o}lder's inequality (cf.
\cite{saloff coste}).

From our experience we knew that we could also derive \textquotedblleft all
the Sobolev inequalities" from the rearrangement inequality%
\begin{equation}
f^{\ast\ast}(t)-f^{\ast}(t)\leq c_{n}t^{1/n}\left\vert \nabla f\right\vert
^{\ast\ast}(t). \label{rea}%
\end{equation}
So it was natural to ask what was the connection between (\ref{gagliardoN})
and (\ref{rea})? The key for us was the method of truncation that had been
devised by Maz'ya (cf. \cite{maz}) and Talenti \cite{tal}, combined with one
natural idea that comes from interpolation theory, or more precisely, from the
work of Marcinkiewicz and Alberto Calder\'{o}n \cite{cal}. The idea simply put
is to do the (smooth) cut-offs using the rearrangement of the function to be
truncated evaluated at a fixed point $t.$ This leads to pointwise
rearrangement inequalities. In this way we could show that (\ref{gagliardoN})
is equivalent to (\ref{rea})! (cf. \cite{mmp})

One advantage of (\ref{rea}) is that since it is a pointwise inequality it
gives all the results, even for non $L^{p}$ norms. This was another innovation
from the work of Alberto Calder\'{o}n \cite{cal}: While previously to prove
interpolation theorems one had to specify the spaces in advance this was not
longer necessary in A. P. Calder\'{o}n's theory. This was crucially important
in our program since the norms than one needs to use to measure the
integrability gains depend on the geometry. We should also point out that for
Euclidean geometries pointwise rearrangement inequalities had been devised and
applied by Talenti (cf. \cite{tal} and the references therein), and since then
we have been applied by many authors working on non-linear PDEs...(cf.
\cite{carlo}, \cite{ra})) In connection with compactness of Sobolev embeddings
connected with the developments in this section we refer to \cite{pu1},
\cite{mmp} and the many references therein.

\section{The Gaussian Inequality: Ledoux's inequality}

The next item on our agenda was a test case: The Gaussian world. The questions
here were: How to formulate the basic inequality (\ref{rea}) in the Gaussian
context? What type of symmetrization was needed to replace the symmetric
rearrangement? What inequality would take the role of the Gagliardo-Nirenberg
inequality in the Gaussian world? The answer to the last question had been
already provided by Ledoux (cf. \cite{led}). Ledoux's inequality is connected
with the use of \ the Maz'ya isoperimetric profile, usually referred to as
isoperimetric profile\footnote{The isoperimetric profile was introduced by
Maz'ya in the sixties and further developed by him in a number of publications
(cf. \cite{maz} and the references therein). Independently, this useful tool
was developed in parallel by geometers (cf. \cite{bayle} and the references
therein) and probabilists (cf. \cite{ledouxbk} and the references therein),
and as we shall see plays an important role in our work formulating Sobolev
pointwise inequalities on rearrangements.} $I_{\gamma},$ associated with
Gaussian measure in $\mathbb{R}^{n},$%
\[
I_{\gamma}(t)=\inf\{\gamma^{+}(A):\gamma(A)=t\},
\]
where $\gamma^{+}$ is the Minkowski content associated with Gaussian measure
defined for Borel sets $A$ by%
\[
\gamma^{+}(A)=\lim_{h\rightarrow0}\frac{\gamma((A_{h}))-\gamma(A)}{h},
\]
where $A_{h}=\{x:d(x,A)<h\},$ and $d$ is the usual $\mathbb{R}^{n}$ distance
between $\{x\}$ and the set $A.$ For Lip functions $f,$ combining the
isoperimetric inequality%
\[
I_{\gamma}\left(  \gamma(\{\left\vert f\right\vert >t\})\right)  \leq
\gamma^{+}(\{\left\vert f\right\vert >t\}),
\]
with the Gaussian co-area formula\footnote{In the general metric case it
becomes the co-area inequality (cf. \cite{bohou}).}, we have
\begin{align}
\int_{0}^{\infty}I_{\gamma}\left(  \gamma(\{\left\vert f\right\vert
>t\})\right)  dt &  \leq\int_{0}^{\infty}\gamma^{+}(\{\left\vert f\right\vert
>t\})dt\label{ledoux}\\
&  =\left\Vert \left\vert \nabla f\right\vert \right\Vert _{L^{1}%
(\mathbb{R}^{n},\gamma)}.\nonumber
\end{align}

This is exactly the same mechanism one can use in the Euclidean
world\footnote{apparently this is true in the whole universe..}. Now, in the
Gaussian world, the isoperimetric profile $I_{\gamma}(t)$ has the following
estimate (with constants independent of the dimension)%
\[
I_{\gamma}(t)\simeq t(\log\frac{1}{t})^{1/2},\, t\in(0,1/2),
\]
while for the Euclidean isoperimetric profile we have%
\[
I_{n}(t)=c_{n}t^{1-1/n},t>0.
\]
Note that in the Euclidean case,
\begin{align*}
\int_{0}^{\infty}I_{n}\left(  \left|  \{\left|  f\right|  >t\}\right|
\right)  dt  &  =\int_{0}^{\infty}I_{n}\left(  \lambda_{f}(t)\right)  dt\\
&  =\int_{0}^{\infty}I_{n}\left(  t\right)  df^{\ast}(t)\\
&  =\int_{0}^{\infty}f^{\ast}(t)dI_{n}\left(  t\right) \\
&  =c_{n}\frac{1}{n^{\prime}}\int_{0}^{\infty}f^{\ast}(t)t^{1/n^{\prime}}%
\frac{dt}{t}.
\end{align*}
In other words, the logarithmic Sobolev inequality\footnote{In this case we
have%
\begin{align*}
\int_{0}^{\infty}I\left(  \gamma(\{\left|  f\right|  >t\})\right)  dt  &
\succeq\int_{0}^{1/2}f^{\ast}(t)d\left(  t\left(  \log\frac{1}{t}\right)
^{1/2}\right) \\
&  \simeq\int_{0}^{1/2}f^{\ast}(t)\left(  \log\frac{1}{t}\right)  ^{1/2}dt.
\end{align*}
} of Ledoux (\ref{ledoux}) is the analogue\footnote{Here again we have to
allow for *generalized* Lorentz spaces since the Gaussian profile although
concave is not increasing. Indeed, $I(t)$ is symmetric about $1/2.$ Also note
that the inequality%
\[
\int_{0}^{\infty}I\left(  \gamma(\{\left|  f\right|  >t\})\right)
dt\leq\left\|  \nabla f\right\|  _{L^{1}(R^{n},\gamma)}%
\]
holds for functions $f$ that do not vanish at the boundary. For example, for
$f=1,$ the right hand is zero and%
\begin{align*}
\int_{0}^{\infty}I\left(  \gamma(\{\left|  f\right|  >t\})\right)  dt  &
=\int_{0}^{1}I\left(  1)\right)  dt\\
&  =\int_{0}^{1}0dt\\
&  =0.
\end{align*}
\par
{}} of the Gagliardo-Nirenberg in Gaussian world. Using the method of
``symmetrization by truncation'' we could then derive the Gaussian version of
(\ref{rea})
\begin{equation}
f_{\gamma}^{\ast\ast}(t)-f_{\gamma}^{\ast}(t)\leq\frac{t}{I_{\gamma}%
(t)}\left|  \nabla f\right|  _{\gamma}^{\ast\ast}(t). \label{rea1}%
\end{equation}
The remarkable fact is that in the Euclidean world (i.e. taking rearrangements
using Lebesgue measure and using the corresponding formula for the Euclidean
isoperimetric profile) this inequality is exactly (\ref{rea})! In fact, these
inequalities are equivalent to the corresponding isoperimetric inequalities in
each of these geometries!

Also note that, the isoperimetric profile automatically selects the spaces
that need to be involved! Our approach in \cite{mamijfa} is based on these
ideas. For example, in the Gaussian case, from (\ref{rea1}) we obtain directly
the following version of a logarithmic Sobolev inequality%
\[
\int\left(  f_{\gamma}^{\ast\ast}(t)-f_{\gamma}^{\ast}(t)\right)  ^{2}\left(
\log\frac{1}{t}\right)  dt\leq\int(\left\vert \nabla f\right\vert _{\gamma
}^{\ast\ast}(t))^{2}dt.
\]
Moreover, since the Gaussian isoperimetric profile does not depend on the
dimension\footnote{It is well known that the sets that realize the
isoperimetric inequality are always hyperspaces: i.e. all but one of the
variables are free.} the inequalities in this case are dimension free! For a
complete development we refer to \cite{mamijfa} and through this paper we
refer to many other important references. Some references connected with this
section..\cite{bart1}, \cite{bart}, \cite{BCR}, \cite{BCR1} \cite{be2},
\cite{bobzeg}, \cite{Bor1}, \cite{bo}, \cite{Ciapick}, \cite{er}, \cite{er1},
\cite{er2}, \cite{gro}, \cite{holeem}, \cite{le1}....

\section{The Metric Case}

When working on the Gaussian inequalities, we realized early on that, with a
suitable definition of modulus of the gradient\footnote{For all Lipschitz
function $f$ on $\Omega,$the modulus of the gradient is defined by
\[
|\nabla f(x)|=\limsup_{d(x,y)\rightarrow0}\frac{|f(x)-f(y)|}{d(x,y)}.
\]
} $\left\vert \nabla f\right\vert ,$ and having at hand an associated co-area
formula\footnote{The strong connection between the co-area formula and Sobolev
embeddings had already been emphasized by Maz'ya in his pionnering fundamental
work in the early sixties (cf. \cite{maz}).}, we could indeed prove
(\ref{rea1}) in the general setting of metric measure spaces (this was
informally first announced in \cite{mamijfa} and more formally in
\cite{mamicomptes}). Fortunately, all the tools that we need to implement this
insight had already been developed by Bobkov-Houdre \cite{bohou}.

In this generalized setting we work with connected metric probability spaces
$(\Omega,\mu)$ The isoperimetric profile $I_{\mu}=I_{(\Omega,d,\mu)}$ is
defined by
\[
I_{(\Omega,d,\mu)}(t)=\inf_{A}\{\mu^{+}(A):\mu(A)=t\},
\]
where $\mu^{+}(A)$ is the perimeter or Minkowski content of the Borel set
$A\subset X,$ defined by
\[
\mu^{+}(A)=\lim\inf_{h\rightarrow0}\frac{\mu\left(  A_{h}\right)  -\mu\left(
A\right)  }{h},
\]
where $A_{h}=\left\{  x\in\Omega:d(x,A)<h\right\}  .$ We assume that $I$ is
continuous, concave, symmetric about $1/2,$ and zero at zero. Further we
assume\footnote{See \cite{he}. Using an approximation argument developed by E.
Milman \cite[Remark 3.3]{mie2} it is possible to prove the main inequalities
of this paper without this assumption (cf. \cite{mamicon} and the forthcoming
\cite{mamilectures}).} that for each $c\in\mathbb{R}$, and each $f\in
Lip(\Omega),$ $|\nabla f(x)|=0,$ $\mu-a.e.$ on the set $\{x:f(x)=c\}.$ The
associated isoperimetric inequality can be formulated as: for all Borel sets
$A$%
\begin{equation}
I(\mu(A))\leq\mu^{+}(A). \label{iso}%
\end{equation}

\subsection{The Gagliardo-Nirenberg-Bobkov-Houdre inequality}

We have the following result due to Bobkov-Houdre (cf. \cite{bohou}),%
\begin{equation}
\int_{0}^{\infty}I_{\mu}(\mu_{f}(t))dt\leq\left\Vert \left\vert \nabla
f\right\vert \right\Vert _{L{^{1}}(\Omega)},\text{ for all }f\in
Lip(\Omega)\Leftrightarrow\text{isoperimetric inequality.}
\label{gagliardoNBH}%
\end{equation}
The reason the isoperimetric inequality is necessary is that given any Borel
set $A$ we can find a sequence of $Lip$ functions $\{f_{n}\}_{n}$ such that
$\left\Vert \left\vert \nabla f_{n}\right\vert \right\Vert _{L{^{1}}(\Omega
)}\rightarrow\mu^{+}(A),$ while $\int_{0}^{\infty}I_{\mu}(\mu_{f_{n}%
}(t))dt\rightarrow I_{\mu}(\mu(A))$ (cf. \cite{bohou}).

\subsection{Generalized P\'{o}lya-Szeg\"{o}}

One big difference between the general probability metric case and the
Gaussian case is the lack of symmetry. In particular, in the Gaussian world
$(\mathbb{R}^{n},\gamma_{n})$ there is a natural choice for a distinguished
rearrangement that replaces the symmetric rearrangement. Let
\[
\phi_{n}(x)=(2\pi)^{-n/2}e^{-\frac{\left|  x\right|  ^{2}}{2}},x\in
\mathbb{R}^{n},\text{ }\Phi(r)=\int_{-\infty}^{r}\phi_{1}(t)dt,\mathbb{R}%
\in\mathbb{R},
\]
then the Gaussian profile is given by (cf. Borell \cite{bo} and
Sudakov-Tsirelson \cite{sut})
\[
I_{\gamma}(t)=\phi_{1}(\Phi^{-1}(t)),\;t\in\lbrack0,1].
\]
The classical Euclidean spherical decreasing rearrangement is replaced by
\[
f_{\gamma_{n}}^{\circ}(x)=f_{\gamma_{n}}^{\ast}(\Phi(x_{1})),
\]
and we have the Erhard analogue of P\'{o}lya-Szeg\"{o} (cf. \cite{er},
\cite{er2}, \cite{mamijfa})%
\begin{equation}
\left|  \nabla\left(  f_{\gamma_{n}}^{\circ}\right)  )\right|  _{\gamma}%
^{\ast\ast}(t)\leq\left|  \nabla f\right|  _{\gamma}^{\ast\ast}(t).
\label{polzgG}%
\end{equation}
In general there is no apparent symmetry and thus no distinguished
rearrangement. This led us to formulate the following alternative inequality
which extends (\ref{polzgG}) to the probability metric case the
P\'{o}lya-Szeg\"{o} inequality in the Euclidean geometry and the Erhard
inequality in the Gaussian case: for all $f\in Lip,$ we have (cf.
\cite{mamiadv})%
\begin{equation}
\int_{0}^{t}\left(  I_{\mu}(\cdot)\frac{d}{dt}(-f_{\mu}^{\ast}(\cdot))\right)
^{\ast}(s)ds\leq\int_{0}^{t}\left|  \nabla f\right|  _{\mu}^{\ast}(s)ds,
\label{polzgGG}%
\end{equation}
where the rearrangement inside the integral on the left hand side is with
respect to Lebesgue measure. The usual formulation of P\'{o}lya-Szeg\"{o} as a
norm inequality follows directly from (\ref{polzgGG}) and the
Calder\'{o}n-Hardy-Littlewood principle. To see this result in detail let us
recall that a rearrangement invariant space $X(\Omega,\mu),$ has a
representation $\bar{X}(0,1)$ such that%
\[
\left\|  f\right\|  _{X(\Omega)}=\left\|  f_{\mu}^{\ast}\right\|  _{X(0,1)}.
\]
Now, since $\left|  \nabla f\right|  _{\mu}^{\ast}=\left(  \left|  \nabla
f\right|  _{\mu}^{\ast}\right)  ^{\ast},$ we see that if (\ref{polzgGG}) holds
then by the Calder\'{o}n-Hardy-Littlewood principle,
\[
\left\|  I_{\mu}(\cdot)\frac{d}{dt}(-f_{\mu}^{\ast}(\cdot))\right\|
_{X(0,1)}\leq\left\|  \left|  \nabla f\right|  _{\mu}^{\ast}\right\|
_{X(0,1)}=\left\|  \left|  \nabla f\right|  _{\mu}^{\ast}\right\|
_{X(\Omega)},
\]
which by abuse of notation (since no confusion can arise) we write as%

\[
\left\|  I_{\mu}(t)\frac{d}{dt}(-f_{\mu}^{\ast}(t))\right\|  _{X}\leq\left\|
\left|  \nabla f\right|  \right\|  _{X}.
\]

The proof of these inequalities follows by smooth truncation. Indeed, if we
apply (\ref{gagliardoNBH}) to the smooth truncations,%
\[
f_{t_{1}}^{t_{2}}(x)=\left\{
\begin{array}
[c]{cc}%
t_{2}-t_{1} & \text{ \ if }\left|  f(x)\right|  \geq t_{2}\\
\left|  f(x)\right|  -t_{1} & \text{ \ if }t_{1}<\left|  f(x)\right|  \leq
t_{2}\\
0 & \text{ \ if }\left|  f(x)\right|  \leq t_{1}%
\end{array}
\right.  ,
\]
and we use the fact that%
\[
\int_{t_{1}}^{t_{2}}I_{\mu}(\mu_{f}(s))ds=\int_{0}^{t_{2}-t_{1}}I_{\mu}%
(\mu_{f_{t_{1}}^{t_{2}}}(s))ds,
\]
then we find that
\[
\int_{t_{1}}^{t_{2}}I_{\mu}(\mu_{f}(s))ds\leq\int_{\{t_{1}<\left|  f\right|
<t_{2}\}}\left|  \nabla f\right|  d\mu.
\]
A careful\footnote{Under the assumption that $f$ is $Lip$ the calculations
below can be justified (cf. \cite{mamiarxiv} \cite{mamiadv}, \cite{leoni})}
``change of variables argument $s\mapsto f_{\mu}^{\ast}(u)$'' (cf.
\cite{mamiadv}) on the left hand side gives%
\[
\int_{t_{1}}^{t_{2}}I_{\mu}(u)\frac{d}{dt}(-f_{\mu}^{\ast})(u)du=\int_{f_{\mu
}^{\ast}(t_{2})}^{f_{\mu}^{\ast}(t_{1})}I_{\mu}(\mu_{f}(s))ds\leq
\int_{\{f_{\mu}^{\ast}(t_{2})<\left|  f\right|  <f_{\mu}^{\ast}(t_{1}%
)\}}\left|  \nabla f\right|  d\mu.
\]
Thus, for any set $E=%
%TCIMACRO{\dbigcup }%
%BeginExpansion
{\displaystyle\bigcup}
%EndExpansion
(a_{i},b_{i})$ union of disjoint intervals, with $\left|  E\right|  =t,$ we
have%
\begin{align*}
\int_{E}I(u)\frac{d}{dt}(-f_{\mu}^{\ast})(u)du  &  \leq%
%TCIMACRO{\dsum \limits_{i}}%
%BeginExpansion
{\displaystyle\sum\limits_{i}}
%EndExpansion
\int_{\{f_{\mu}^{\ast}(b_{i})<\left|  f\right|  <f_{\mu}^{\ast}(a_{i}%
)\}}\left|  \nabla f\right|  d\mu\\
&  =\int_{%
%TCIMACRO{\dbigcup \limits_{i}}%
%BeginExpansion
{\displaystyle\bigcup\limits_{i}}
%EndExpansion
\{f_{\mu}^{\ast}(b_{i})<\left|  f\right|  <f_{\mu}^{\ast}(a_{i})\}}\left|
\nabla f\right|  (s)d\mu\\
&  \leq\int_{0}^{%
%TCIMACRO{\dsum }%
%BeginExpansion
{\displaystyle\sum}
%EndExpansion
\mu\{f_{\mu}^{\ast}(b_{i})<\left|  f\right|  <f_{\mu}^{\ast}(a_{i})\}}\left|
\nabla f\right|  _{\mu}^{\ast}(s)d\mu\\
&  \leq\int_{0}^{%
%TCIMACRO{\dsum }%
%BeginExpansion
{\displaystyle\sum}
%EndExpansion
(b_{i}-a_{i})}\left|  \nabla f\right|  _{\mu}^{\ast}(s)d\mu\\
&  =\int_{0}^{t}\left|  \nabla f\right|  _{\mu}^{\ast}(s)d\mu.
\end{align*}
By approximation,%
\begin{align*}
\int_{0}^{t}(I(\cdot)\frac{d}{dt}(-f^{\ast})(\cdot))^{\ast}(u)du  &
=\sup_{\left|  E\right|  =t}\int_{E}I(u)\frac{d}{dt}(-f^{\ast})(u)du\\
&  \leq\int_{0}^{t}\left|  \nabla f\right|  _{\mu}^{\ast}(s)d\mu,
\end{align*}
and (\ref{polzgGG}) follows.

\subsection{Maz'ya-Talenti inequality}

Likewise, with a similar direct proof, we can extend a well known inequality
apparently obtained first by Maz'ya \cite{maz1} and independently, but later,
by Talenti (cf. \cite{mmp} and the references therein). In our setting this
inequality takes the following form (cf. \cite{mamiadv})%
\begin{equation}
I_{\mu}(t)\frac{d}{dt}(-f_{\mu}^{\ast}(t))\leq\frac{d}{dt}\int_{\{\left\vert
f\right\vert >f^{\ast}(t)\}}\left\vert \nabla(f)\right\vert _{\mu}^{\ast
}(s)ds. \label{maztal}%
\end{equation}
In fact, it is also easy to see that the argument given above shows that
(\ref{maztal}) implies (\ref{polzgGG}). Indeed, starting with (\ref{maztal}),
and considering first sets which are unions of disjoint intervals, we arrive
at
\[
\sup_{\left\vert E\right\vert =t}\int_{E}I_{\mu}(u)\frac{d}{dt}(-f_{\mu}%
^{\ast})(u)du\leq\int_{0}^{t}\left\vert \nabla f\right\vert _{\mu}^{\ast
}(s)d\mu.
\]

\subsection{Basic rearrangement inequality}

Next we observe that the previous arguments with intervals also show (cf.
\cite{mamiadv}) that, integrating by parts, we obtain (recall that under our
assumptions $f_{\mu}^{\ast}$ is absolutely continuous (cf. \cite{leoni})%
\[
f_{\mu}^{\ast\ast}(t)-f_{\mu}^{\ast}(t)=\frac{1}{t}\int_{0}^{t}s\frac{d}%
{dt}(-f_{\mu}^{\ast})(s)ds.
\]
Consequently,%
\begin{equation}
f_{\mu}^{\ast\ast}(t)-f_{\mu}^{\ast}(t)\leq\frac{1}{t}\frac{t}{I_{\mu}(t)}%
\int_{0}^{t}I_{\mu}(s)\frac{d}{dt}(-f_{\mu}^{\ast})(s)ds\text{ \ \ (since
}\frac{s}{I_{\mu}(s)}\text{ is increasing).} \label{laotra}%
\end{equation}
By now, we have discussed several inequalities that can be used to control
$\int_{0}^{t}I_{\mu}(s)(-f_{\mu}^{\ast}(s))^{\prime}ds.$ For example, using
the definition of rearrangement and the generalized P\'{o}lya-Szeg\"{o}
inequality (\ref{polzgGG}), we have%
\begin{align*}
\int_{0}^{t}I_{\mu}(s)\frac{d}{dt}(-f_{\mu}^{\ast})(s)ds  &  \leq\int_{0}%
^{t}\left(  I_{\mu}(\cdot)\frac{d}{dt}(-f_{\mu}^{\ast})(\cdot)\right)  ^{\ast
}(s)ds\\
&  \leq\int_{0}^{t}\left|  \nabla f\right|  _{\mu}^{\ast}(s)ds.
\end{align*}
Consequently, combining the last inequality with (\ref{laotra}) we arrive to
the familiar%
\begin{align}
f_{\mu}^{\ast\ast}(t)-f_{\mu}^{\ast}(t)  &  \leq\frac{1}{I_{\mu}(t)}\int
_{0}^{t}\left|  \nabla f\right|  _{\mu}^{\ast}(s)ds\nonumber\\
&  =\frac{t}{I(t)}\left|  \nabla f\right|  _{\mu}^{\ast\ast}(t).
\label{metricas}%
\end{align}
Using a standard approximation argument (like the one outlined above) (cf.
\cite{mamiadv}, \cite{mamicon}) we see that (\ref{metricas}) implies the
isoperimetric inequality (\ref{iso}).

Let us also remark that is fairly easy to give a direct proof of
(\ref{metricas}) (cf. Section \ref{seccoulhon} below)

\subsection{Self Improvement}

Under stronger assumptions on the isoperimetric profile we can get stronger
inequalities. For example, suppose that
\[
\int_{t}^{1}\frac{I_{\mu}(s)}{s}\frac{ds}{s}\leq c\frac{I_{\mu}(t)}%
{t},\;\;t\in(0,1).
\]
Then, for all $f\in Lip(\Omega),$
\begin{equation}
\int_{0}^{t}(f_{\mu}^{\ast\ast}(s)-f_{\mu}^{\ast}(s))\frac{I(s)}{s}ds\leq
\int_{0}^{t}\left(  \frac{I(\cdot)}{\left(  \cdot\right)  }[f_{\mu}^{\ast\ast
}(\cdot)-f_{\mu}^{\ast}(\cdot)]\right)  ^{\ast}ds\leq C\int_{0}^{t}\left|
\nabla f\right|  _{\mu}^{\ast}(s)ds. \label{l1}%
\end{equation}
This extends the sharp form of the Gagliardo-Nirenberg inequality to metric
spaces that have ``Euclidean like profiles''. Indeed, suppose that
$\frac{I(s)}{s}\simeq s^{1/n^{\prime}},$ then for $f\in Lip_{0}(\Omega),$ we
have%
\[
\int_{0}^{t}(f_{\mu}^{\ast\ast}(s)-f_{\mu}^{\ast}(s))s^{1/n^{\prime}}\frac
{ds}{s}\leq c_{n}\int_{0}^{t}\left|  \nabla f\right|  _{\mu}^{\ast}(s)ds.
\]
Consequently, letting $t\rightarrow\infty,$%
\[
\int_{0}^{\infty}(f_{\mu}^{\ast\ast}(s)-f_{\mu}^{\ast}(s))s^{1/n^{\prime}%
}\frac{ds}{s}\leq c_{n}\left\|  \nabla f\right\|  _{1}.
\]

\subsection{Poincar\'{e} inequalities and generalized P\'{o}lya-Szeg\"{o}}

Once again we follow \cite{mamiadv}. To discuss Poincar\'{e} inequalities in
the metric setting we introduced the Hardy isoperimetric operators (cf.
\cite{mamiadv}). Let us define%
\[
Q_{I}(f)(t)=\chi_{(0,1/2)}(t)\int_{t}^{1/2}f(s)\frac{ds}{I_{\mu}(s)}.
\]
Suppose that for positive functions supported on $(0,1/2)$%
\[
\left\Vert Q_{I}(f)\right\Vert _{X}\leq\left\Vert Q_{I}\right\Vert
_{X\rightarrow X}\left\Vert f\right\Vert _{X}.
\]
Then, the following Poincar\'{e} inequality holds%
\[
\left\Vert g-\int_{\Omega}g\right\Vert _{X}\preceq\left\Vert Q_{I}\right\Vert
_{X\rightarrow X}\left\Vert \left\vert \nabla g\right\vert \right\Vert _{X}.
\]

Indeed, for $g\in Lip(\Omega),$ $t\in(0,1/2)$ we can write$,$%
\begin{align}
g_{\mu}^{\ast}(t)-g_{\mu}^{\ast}(1/2)  &  =\int_{t}^{1/2}\left(  -g_{\mu
}^{\ast}(s)\right)  ^{\prime}ds\nonumber\\
&  =\int_{t}^{1/2}I_{\mu}(s)\left(  -g_{\mu}^{\ast}(s)\right)  ^{\prime}%
\frac{ds}{I_{\mu}(s)}\nonumber\\
&  =Q_{I}(I_{\mu}(.)\left(  -g_{\mu}^{\ast}(.)\right)  ^{\prime}).
\label{robusta}%
\end{align}
Thus,%
\begin{equation}
g_{\mu}^{\ast}(t)=Q_{I}(I_{\mu}(.)\left(  -g_{\mu}^{\ast}(.)\right)  ^{\prime
})+g_{\mu}^{\ast}(1/2) \label{from}%
\end{equation}
Now,%
\[
\left\|  g_{\mu}^{\ast}\right\|  _{X}\leq\left\|  g_{\mu}^{\ast}\chi
_{(0,1/2)}\right\|  _{X}+\left\|  g_{\mu}^{\ast}\chi_{(1/2,1)}\right\|  _{X},
\]
and
\begin{align*}
\left\|  g_{\mu}^{\ast}\chi_{(1/2,1)}\right\|  _{X}  &  \leq g_{\mu}^{\ast
}(1/2)\left\|  \chi_{(1/2,1)}\right\|  _{X}\\
&  =g_{\mu}^{\ast}(1/2)\left\|  \chi_{(0,1/2)}\right\|  _{X}\\
&  \leq\left\|  g_{\mu}^{\ast}\chi_{(0,1/2)}\right\|  _{X},
\end{align*}
yields%
\[
\left\|  g_{\mu}^{\ast}\right\|  _{X}\preceq\left\|  g_{\mu}^{\ast}%
\chi_{(0,1/2)}\right\|  _{X}%
\]
which combined with (\ref{from}) gives%
\begin{align*}
\left\|  g_{\mu}^{\ast}\right\|  _{X}  &  \preceq\left\|  g_{\mu}^{\ast}%
\chi_{(0,1/2)}\right\|  _{X}\\
&  \leq\left\|  Q_{I}(I_{\mu}(\cdot)\frac{d}{dt}(-g_{\mu}^{\ast}%
)(\cdot)\right\|  _{X}+\left\|  g_{\mu}^{\ast}(1/2)\right\|  _{X}\\
&  \leq\left\|  Q_{I}\right\|  _{X\rightarrow X}\left\|  (I_{\mu}(\cdot
)\frac{d}{dt}(-g_{\mu}^{\ast}))\right\|  _{X}+\left\|  g_{\mu}^{\ast
}(1/2)\right\|  _{X}\\
&  \leq\left\|  Q_{I}\right\|  _{X\rightarrow X}\left\|  \left|  \nabla
g\right|  \right\|  _{X}+g_{\mu}^{\ast}(1/2)\left\|  1\right\|  _{X}.
\end{align*}
Now, by Chebyshev's inequality $g_{\mu}^{\ast}(1/2)\leq2\left\|  g\right\|
_{1}.$ Therefore,%
\begin{align*}
\left\|  g\right\|  _{X}  &  =\left\|  g_{\mu}^{\ast}\right\|  _{X}\\
&  \preceq\left\|  Q_{I}\right\|  _{X\rightarrow X}(\left\|  \left|  \nabla
g\right|  \right\|  _{X}+\left\|  g\right\|  _{1}).
\end{align*}
Applying the previous inequality to $g-\int_{\Omega}g,$%
\[
\left\|  g-\int_{\Omega}g\right\|  _{X}\preceq\left\|  Q_{I}\right\|
_{X\rightarrow X}(\left\|  \left|  \nabla g\right|  \right\|  _{X}+\left\|
g-\int_{\Omega}g\right\|  _{1}),
\]
and combining with Cheeger's inequality yields%
\begin{align*}
\left\|  g-\int_{\Omega}g\right\|  _{1}  &  \preceq\left\|  \nabla g\right\|
_{1}\\
&  \leq\left\|  \nabla g\right\|  _{X}.
\end{align*}
Thus we finally arrive at%
\[
\left\|  g-\int_{\Omega}g\right\|  _{X}\preceq\left\|  Q_{I}\right\|
_{X\rightarrow X}\left\|  \left|  \nabla g\right|  \right\|  _{X}.
\]

\subsection{Identity associated with Poincar\'{e}}

For the sake of completeness we note that the previous discussion also shows
that for $g\in Lip(\Omega)$ we have (cf. \ref{robusta})%
\[
g_{\mu}^{\ast}(t)-g_{\mu}^{\ast}(1/2)=Q_{I}(I_{\mu}(\cdot)\frac{d}{dt}%
(-g_{\mu}^{\ast})(\cdot)).
\]

This section is based entirely on \cite{mamiadv} but builds on a long list of
important contributors we refer to \cite{mamiadv} for a more complete list of
references. For versions of P\'{o}lya-Szeg\"{o} in the classical format (norm
inequality) we refer to the books \cite{kesavan}, \cite{leoni} as well as
Talenti \cite{tal}, Almgren-Lieb \cite{alm}, and the references therein. We
should also mention the work of the French geometer Gallot \cite{gallot} on
symmetrization inequalities on manifolds, which as we discovered a posteriori,
is close in spirit to our development here. The important paper by D. Bakry,
T. Coulhon, M. Ledoux and L. Saloff-Coste \cite{bakr} opened our eyes early on
to the possibilities afforded by the masterful use of cut-offs.

\section{Coulhon Inequalities\label{seccoulhon}}

The previous section dealt with Sobolev inequalities in probability metric
spaces. The next natural questions are: how to deal with more general metric
measure spaces? Secondly: what is the rearrangement version of $L^{p}$ Sobolev
inequalities for $p>1?$ In this situation we need to replace the isoperimetric
profile by a suitable capacitary profile and the formulation of the
inequalities is somewhat more complicated (cf. \cite{mamicon}). However, one
can take a different tack (cf. \cite{mamicoulhon}) and work with an equivalent
formulation of the Sobolev inequalities due to Coulhon (cf. \cite{cou1},
\cite{cou2}, \cite{cou}) and Bakry-Coulhon-Ledoux \cite{bakr}.

Let $(\Omega,d,\mu)$ be a metric measure space, let $p\in\lbrack1,\infty],$
and let $\phi$ be an increasing function on the positive half line. We
consider Coulhon inequalities of the form
\[
(S_{\phi}^{p})\;\;\;\;\left\Vert f\right\Vert _{p}\leq\phi(\left\Vert
f\right\Vert _{0})\left\Vert \left\vert \nabla f\right\vert \right\Vert
_{p},\text{ }f\in Lip_{0}(\Omega),
\]
where $\left\Vert f\right\Vert _{0}=\mu\{$\textrm{support} $(f)\}.$

An important point to keep in mind is that the function $\phi$ must then be
connected to the geometry; but now we do not specify what this connection is
in advance! In \cite{mamicon} we find rearrangement inequalities that
characterize $(S_{\phi}^{p}).$ As was to be expected the rearrangement
inequalities incorporate in their formulation the function $\phi$.

To see the connection with the previous discussions let us consider the case
$p=1.$ We shall show that this case is connected with the isoperimetric
inequality. Before doing so, let us mention that our approach in this section
is independent and, indeed, can be seen as an alternative route to
rearrangement inequalities discussed in previous sections. Let us then
consider the connection between $(S_{\phi}^{1})$ and a rearrangement
inequality of the form%
\begin{equation}
f_{\mu}^{\ast\ast}(t)-f_{\mu}^{\ast}(t)\leq\phi(t)\left|  \nabla f\right|
_{\mu}^{\ast\ast}(t),\text{ }f\in Lip_{0}(\Omega). \label{verde1}%
\end{equation}
Suppose (\ref{verde1}) holds where $\phi$ is a given increasing continuous
function. Let $t>0;$ multiplying both sides of (\ref{verde1}) by $t,$ we
obtain
\[
t\left(  f_{\mu}^{\ast\ast}(t)-f_{\mu}^{\ast}(t)\right)  \leq\phi(t)\int
_{0}^{t}\left|  \nabla f\right|  _{\mu}^{\ast}(s)ds.
\]
Since formally $f_{\mu}^{\ast}(t)=\mu_{f}^{-1}(t),$ drawing a diagram it is
easy to see that
\begin{align*}
t\left(  f_{\mu}^{\ast\ast}(t)-f_{\mu}^{\ast}(t)\right)   &  =\int_{0}%
^{t}f_{\mu}^{\ast}(s)ds-tf_{\mu}^{\ast}(t)\\
&  =\int_{f^{\ast}(t)}^{\infty}\mu_{f}(s)ds.
\end{align*}
Consequently, if we let $t=\left\|  f\right\|  _{0},$ we see that $f_{\mu
}^{\ast}(\left\|  f\right\|  _{0})=0,$ $\int_{f_{\mu}^{\ast}(\left\|
f\right\|  _{0})}^{\infty}\mu_{f}(s)ds=\left\|  f\right\|  _{1},$ and
$\int_{0}^{\left\|  f\right\|  _{0}}\left|  \nabla f\right|  _{\mu}^{\ast
}(s)ds=\left\|  \left|  \nabla f\right|  \right\|  _{1}.$ Thus,
\begin{equation}
\left\|  f\right\|  _{1}\leq\phi(\left\|  f\right\|  _{0})\left\|  \left|
\nabla f\right|  \right\|  _{1}. \label{espada}%
\end{equation}
Suppose on the other hand that an $(S_{\tilde{\phi}}^{1})$ condition holds.
For $f\in Lip_{0}(\Omega),$ and for $t>0,$ let us apply the $(S_{\phi}^{1})$
condition to $\left[  f-f_{\mu}^{\ast}(t)\right]  _{+}.$ We compute,%
\[
\left\|  \left[  f-f_{\mu}^{\ast}(t)\right]  _{+}\right\|  _{0}=\mu\{f>f_{\mu
}^{\ast}(t)\}\leq t,
\]
moreover, since $\int_{\{f=f^{\ast}(t)\}}\left|  \nabla\left[  f(x)-f_{\mu
}^{\ast}(t)\right]  \right|  dx=0,$ we have
\[
\left\|  \nabla\left[  f-f_{\mu}^{\ast}(t)\right]  _{+}\right\|  _{L^{1}}%
=\int_{\{f>f^{\ast}(t)\}}\left|  \nabla f\right|  _{\mu}^{\ast}(s)ds.
\]
We also have%
\begin{align*}
\left\|  \left[  f-f_{\mu}^{\ast}(t)\right]  _{+}\right\|  _{1}  &
=\int_{\{f>f^{\ast}(t)\}}\left[  f(s)-f_{\mu}^{\ast}(t)\right]  _{+}d\mu(s)\\
&  =\int_{0}^{\infty}[f_{\mu}^{\ast}(x)-f_{\mu}^{\ast}(t)]_{+}\,dx\\
&  =\int_{0}^{t}(f_{\mu}^{\ast}(x)-f_{\mu}^{\ast}(t))\,dx\\
&  =t(f_{\mu}^{\ast\ast}(t)-f_{\mu}^{\ast}(t)).
\end{align*}
Inserting these calculations in (\ref{espada}), we find%
\[
t(f_{\mu}^{\ast\ast}(t)-f_{\mu}^{\ast}(t))\leq t\phi(t)\left(  \frac{1}{t}%
\int_{0}^{t}\left|  \nabla f\right|  _{\mu}^{\ast}(s)ds\right)  .
\]
In other words, we have shown that
\[
(S_{\phi}^{1})\Leftrightarrow(f_{\mu}^{\ast\ast}(t)-f_{\mu}^{\ast}(t))\leq
\phi(t)\left|  \nabla f\right|  _{\mu}^{\ast\ast}(t),f\in Lip_{0}(\Omega).
\]
In particular, if $I_{\mu}(t)$ is concave then $\phi(t)=\frac{t}{I_{\mu}(t)}$
is increasing, and we have
\[
(S_{\frac{t}{I_{\mu}(t)}}^{1})\Leftrightarrow(f_{\mu}^{\ast\ast}(t)-f_{\mu
}^{\ast}(t))\leq\frac{t}{I_{\mu}(t)}\left|  \nabla f\right|  _{\mu}^{\ast\ast
}(t),f\in Lip_{0}(\Omega).
\]
Let us now show that the function $\frac{t}{I_{\mu}(t)}$ is optimal. Suppose
that an $(S_{\phi}^{1})$ condition holds, then it is easy to see, by
approximation, that $\frac{t}{I_{\mu}(t)}\leq\phi(t)$ (cf. \cite{mamicoulhon}%
)$.$ The fact that $(S_{\frac{t}{I_{\mu}(t)}}^{1})$ itself holds is a direct
consequence of the co-area inequality for $Lip$ functions on metric spaces
(cf. \cite{bohou}). It is instructive to see the details. Suppose that $f\in
Lip_{0}(\Omega),$ then from $I_{\mu}(\mu(A))\leq\mu^{+}(A),$ for all Borel
sets, we have the Bobkov-Houdre inequality%
\[
\int_{0}^{\infty}I_{\mu}(\mu_{f}(t))dt\leq\left\|  \left|  \nabla f\right|
\right\|  _{1}.
\]
Now, since $\frac{I_{\mu}(t)}{t}$ decreases, and we obviously have $\mu
_{f}(t)\leq\left\|  f\right\|  _{0},$ we see that
\begin{align*}
\int_{0}^{\infty}I_{\mu}(\mu_{f}(t))dt  &  =\int_{0}^{\infty}\mu_{f}%
(t)\frac{I_{\mu}(\mu_{f}(t))}{\mu_{f}(t)}dt\\
&  \geq\int_{0}^{\infty}\mu_{f}(t)dt\frac{I_{\mu}(\left\|  f\right\|  _{0}%
)}{\left\|  f\right\|  _{0}}.
\end{align*}
Combining these inequalities, we see that for the choice of $\phi(t)=\frac
{t}{I_{\mu}(t)}$ we do indeed have%
\[
\left\|  f\right\|  _{1}=\int_{0}^{\infty}\mu_{f}(t)dt\leq\phi(\left\|
f\right\|  _{0})\left\|  \left|  \nabla f\right|  \right\|  _{1}.
\]

For example, in the Euclidean space $\mathbb{R}^{n}$, $I(t)=d_{n}t^{1-1/n},$
$\phi(t)\simeq t^{1/n}$ and the best possible $(S_{\phi}^{1})$ inequality can
be written as
\[
\left\|  f\right\|  _{1}\leq c_{n}\left\|  f\right\|  _{0}^{1/n}\left\|
\left|  \nabla f\right|  \right\|  _{1}.
\]
With the choice of $\phi(t)=t^{1/n}$, and $p=2,$ Coulhon \cite{cou1} shows the
equivalence of the $S_{\phi}^{2}$ inequality with the classical Nash
inequality
\[
\left\|  f\right\|  _{2}^{1+2/n}\leq c_{n}\left\|  f\right\|  _{1}%
^{2/n}\left\|  \left|  \nabla f\right|  \right\|  _{2}.
\]

The general characterization of Coulhon inequalities in terms of
rearrangements is given by the following (cf. \cite{mamicoulhon})

\begin{theorem}
\label{teo}Let $(\Omega,d,\mu)$ be a connected Borel metric measure space as
described above, and let $p\in\lbrack1,\infty).$ The following statements are equivalent

\begin{enumerate}
\item $(S_{\phi}^{p})$ holds, i.e.
\[
\left\Vert f\right\Vert _{p}\leq\phi(\left\Vert f\right\Vert _{0})\left\Vert
\left\vert \nabla f\right\vert \right\Vert _{p},\text{ for all }f\in
Lip_{0}(\Omega).
\]

\item Let $k\in\mathbb{N}$ be such that $k<p\leq k+1,$ then for all $f\in
Lip_{0}(\Omega)$
\begin{equation}
\left(  \frac{f_{\left(  p\right)  }^{\ast\ast}(t)}{\phi_{(p)}(t)}\right)
^{1/p}-\left(  \frac{f_{(p)}^{\ast}(t)}{\phi_{(p)}(t)}\right)  ^{1/p}%
\leq2^{\frac{k+1}{p}-1}\left(  \left\vert \nabla f\right\vert _{(p)}^{\ast
\ast}(t)\right)  ^{1/p}, \label{norma}%
\end{equation}
where
\[
f_{(p)}^{\ast}(t)=\left(  f^{\ast}(t)\right)  ^{p},\text{ }f_{(p)}^{\ast\ast
}(t)=\frac{1}{t}\int_{0}^{t}f_{(p)}^{\ast}(s)ds,\text{ }\phi_{(p)}(t)=\left(
\phi(t)\right)  ^{p}.
\]

\end{enumerate}
\end{theorem}

We note that for $p=1$ the inequality (\ref{norma}) of Theorem \ref{teo}
coincides with (\ref{metricas}). This new characterization for $p\geq1$ is
independent of \cite{mamiadv}, and, in fact, the proof we gave above provides
a new approach to (\ref{metricas}) as well.

For the details of the proof of Theorem \ref{teo} we must refer to
\cite{mamicoulhon}.

We cannot resist to make a connection between Coulhon's inequalities and the
theory of factorization of operators. Local operators (cf. \cite{pis},
\cite{allen}, \cite{milocal} and the references therein) satisfy conditions of
the form%
\[
\left\Vert Tf\right\Vert _{Y}\leq\phi(\left\Vert f\right\Vert _{0})\left\Vert
f\right\Vert _{X}.
\]
Self-improvements in this setting are expressed via factorization and change
of density. We think that the factorization of Sobolev inequalities could be
an interesting line of investigation to pursue .

\section{Connection with the work of Emanuel Milman}

In this section, following \cite{mamiadv}, we consider a connection with the
work of Emanuel Milman\footnote{We refer to E. Milman's papers for an account
of the history of the problem. Emanuel, who belongs to the Milman family of
mathematicians that includes David (grandfather), Vitali (father), Pierre
(uncle) (cf. \cite{Vmilman}), is no direct relation to Mario Milman.} (cf.
\cite{mie1}, \cite{MiE}, \cite{mie2}, \cite{miemie}, \cite{miemaz},
\cite{milcon})).

For metric measure spaces $(\Omega,d,\mu)$ obtained from a $C^{\infty}$
complete oriented $n-$dimensional Riemannian manifold $(M,g)$, where $d$ is
the induced geodesic distance and $\mu$ is absolutely continuous with respect
to dvol$_{M},$ it is known that the corresponding isoperimetric profile
satisfies (cf. \cite{bayle}) that $I_{(\Omega,d,\mu)}(t)$ is continuous,
$I_{(\Omega,d,\mu)}(t)>0$ for $t\in(0,1),$ and moreover,%
\[
I_{(\Omega,d,\mu)}(t)=I_{(\Omega,d,\mu)}(1-t),\forall t\in\lbrack0,1].
\]
E. Milman further assumes some convexity conditions: $d\mu=e^{-\Psi}dvol_{M},$
where $\Psi$ is such that $\Psi\in C^{2}(M),$ and as tensor fields
Ric$_{g}+Hess_{g}(\Psi)\geq0$ on $M.$ In this case it then follows that
$I_{(\Omega,d,\mu)}$ is also concave (cf. \cite{mie1} and the extensive list
of references therein). Under such conditions, E. Milman shows the equivalence
of Cheeger's inequality, Poincar\'{e}'s inequality and concentration
inequalities! More precisely, using a variety of different tools, including
the semigroup approach of Ledoux, E. Milman has shown that (cf. also Ledoux's
\cite{ledcon} streamlined approach to E. Milman's results in \cite{miemie}).

\begin{theorem}
\label{teoem1}(E. Milman) Let $(\Omega,d,\mu)$ be a metric probability space
satisfying E. Milman's convexity conditions. Then following statements are equivalent

\noindent(E1) Cheeger's inequality: there exists a positive constant $C$ such
that%
\[
I_{(\Omega,d,\mu)}\geq Ct,\ \ \ t\in(0,1/2].
\]
(E2) Poincar\'{e}'s inequality: there exists a positive constant $P$ such that
for all $f\in Lip(\Omega),$%
\[
\left\Vert f-m_{e}\right\Vert _{L^{2}(\Omega)}\leq P\left\Vert \left\vert
\nabla f\right\vert \right\Vert _{L^{2}(\Omega)}.
\]
(E3) Exponential concentration: there exist positive constants $c_{1},c_{2}$
such that for all $f\in Lip(\Omega)$ with $\left\Vert f\right\Vert
_{Lip(\Omega)}\leq1,$%
\[
\mu\{\left\vert f-m_{e}\right\vert >t\}\leq c_{1}e^{-c_{2}t},\text{ \ }%
t\in(0,1).
\]
(E4) First moment inequality: there exists a positive constant $F$ such that
for all $f\in Lip(\Omega)$ with $\left\Vert f\right\Vert _{Lip(\Omega)}\leq1,$%
\[
\left\Vert f-m_{e}\right\Vert _{L^{1}(\Omega)}\leq F.
\]

\end{theorem}

\begin{remark}
For the optimal relationship between the concentration profile and Cheeger's
constant under suitable convexity conditions see \cite{milcon}.
\end{remark}

\subsection{Isoperimetric Hardy type}

We single out a class of metric probability spaces that are suitable for our
analysis (cf. \cite{mamiadv}, \cite{mamimazya}).

\begin{definition}
\label{def:isohar}We shall say that a probability metric space $(\Omega
,d,\mu)$ is of isoperimetric Hardy type if for any given isoperimetric
estimator $I,$ the following are equivalent for all r.i. spaces $X=X(\Omega)$,
$Y=Y(\Omega)$.

\begin{enumerate}
\item There exists a constant $c=c(X,Y)$ such that for all $f\in Lip(\Omega)$
\[
\left\|  f-\int_{\Omega}fd\mu\right\|  _{Y}\leq c\left\|  \left|  \nabla
f\right|  \right\|  _{X}.
\]

\item There exists a constant $c_{1}=c_{1}(X,Y)>0$ such that for all positive
functions $f\in\bar{X},$ with $supp(f)\subset(0,1/2)$ we have
\[
\left\Vert Q_{I}f\right\Vert _{\bar{Y}}\leq c_{1}\left\Vert f\right\Vert
_{\bar{X}},
\]
where $Q_{I}$ is the isoperimetric Hardy operator
\[
Q_{I}f(t)=\chi_{(0,1/2)}(t)\int_{t}^{1/2}f(s)\frac{ds}{I(s)}.
\]

\end{enumerate}
\end{definition}

For spaces of isoperimetric type it is possible to give a very simple proof of
the E. Milman's equivalences (cf. \cite{mamiadv}).

\begin{theorem}
\label{teoema1}Suppose that $(\Omega,d,\mu)$ is a metric probability space of
isoperimetric Hardy type. Then
\[
(E1)\Leftrightarrow(E2)\Leftrightarrow(E3)\Leftrightarrow(E4).
\]

\end{theorem}

\begin{example}
All the model spaces studied in \cite{mamiadv} (including Gaussian space (cf.
\cite{mamijfa}) are of Hardy isoperimetric type.
\end{example}

For further results connecting our work with E. Milman's work we refer to
\cite{mamiadv} and \cite{mamicon}. For example, the following result of E.
Milman can be understood in the context of Hardy isoperimetric type (cf.
\cite{mamiadv}).

\begin{theorem}
\label{teoem2}Let $(\Omega,d,\mu)$ be a space satisfying E. Milman's convexity
conditions. Let $1\leq q<\infty,$ and let $N$ be a Young's function such that
$\frac{N(t)^{1/q}}{t}$ is non-decreasing, and there exists $\alpha>\max
\{\frac{1}{q}-\frac{1}{2},0\}$ such that $\frac{N(t^{\alpha})}{t}$
non-increasing. Then, the following statements are equivalent:

\noindent(E5) $(L_{N},L^{q})$ Poincar\'{e} inequality holds: there exists a
positive constant $P$ such that for all $f\in Lip(\Omega)$%
\[
\left\Vert f-m_{e}\right\Vert _{L_{N}(\Omega)}\leq P\left\Vert \left\vert
\nabla f\right\vert \right\Vert _{L^{q}(\Omega)}.
\]
(E6) Any isoperimetric profile estimator $I$ satisfies: there exists a
constant $c>0$ such that $I(t)\geq c\frac{t^{1-1/q}}{N^{-1}(1/t)},$
\ $t\in(0,1/2].$
\end{theorem}

E. Milman's work can also be seen as providing a program to unify Sobolev
inequalities in different geometries. For more on this we refer to our paper
\cite{mamicon} where, generalizing E. Milman's work, we in particular show why
Lorentz spaces appear as optimal target spaces for Sobolev embeddings.

\section{Transference and dimensionless inequalities}

In this section we follow \cite{mamiadv} and \cite{mamidim} to show how our
pointwise inequalities can be used to transfer Sobolev inequalities from one
geometry to another. The transference is of special interest when it is
implemented to replace Sobolev inequalities that carry dimensional constants
by weaker, but dimensionless, Sobolev inequalities. Inequalities independent
of the dimension play an increasingly important role in approximation theory
and its applications (cf. \cite{griebel}). To see how our pointwise
inequalities are relevant for this task let us recall that our typical Sobolev
inequality on a metric probability space $(\Omega,d,\mu),$ takes the form%
\begin{equation}
\left\|  \left(  f_{\mu}^{\ast\ast}(t)-f_{\mu}^{\ast}(t)\right)  \frac{I_{\mu
}(t)}{t}\right\|  _{\bar{X}}\leq c\left\|  \left|  \nabla f\right|  _{\mu
}^{\ast\ast}\right\|  _{\bar{X}}. \label{latia}%
\end{equation}
Now suppose that the metric probability space $(\Omega,d,\mu)$ is of
``Gaussian isoperimetric type'', that is suppose that for some universal
constant independent of the dimension, it holds
\[
I_{(\Omega,\mu)}(t)\succeq t\left(  \log\frac{1}{t}\right)  ^{\frac{1}{2}%
},\text{ on }\left(  0,\frac{1}{2}\right)  ;
\]
then we can obviously replace $\frac{I_{\mu}(t)}{t}$ by $\left(  \log\frac
{1}{t}\right)  ^{\frac{1}{2}}$ in (\ref{latia}), and in this fashion
*transfer* the Gaussian inequality to $(\Omega,d,\mu):$%
\begin{equation}
\left\|  \left(  f_{\mu}^{\ast\ast}(t)-f_{\mu}^{\ast}(t)\right)  \left(
\log\frac{1}{t}\right)  ^{\frac{1}{2}}\right\|  _{\bar{X}}\leq c\left\|
\left|  \nabla f\right|  _{\mu}^{\ast\ast}\right\|  _{\bar{X}}.
\label{lavieja}%
\end{equation}
This argument shows how the Gaussian log Sobolev inequalities can be
transferred to $(\Omega,d,\mu)$ with constants independent of the dimension.
In particular, since $Q_{n},$ the open unit cube in $\mathbb{R}^{n},$ is of
Gaussian type with constant equal to $1$ (cf. \cite{Ros})$,$ the Gaussian
Sobolev inequalities can be transferred to $Q_{n}$, with constants independent
of the dimension.

This answered a question of Triebel \cite{tri} (cf. \cite{mamiadv},
\cite{krb1}, \cite{krb2}, \cite{mamidim}, \cite{tri1}) and the references
therein). For example, the non-homogeneous form of these results take the
following form%
\begin{equation}
\left\|  f\right\|  _{L^{q}(LogL)^{q/2}(Q_{n})}\leq C(q)\left(  \left\|
\left|  \nabla f\right|  \right\|  _{L^{q}(Q_{n})}+\left\|  f\right\|
_{L^{q}(Q_{n})}\right)  ,\text{ }f\in C_{0}^{\infty}(Q_{n}).
\label{launoprima}%
\end{equation}
In turn, this result was recently improved by Krbec-Fiorenza-Schmeisser
\cite{FKS}, using the spaces $L_{(q,q^{\prime}}(Q_{n}),$ originally introduced
by Iwaniec-Sbordone-Fiorenza (cf. \cite{iw}, \cite{fio1}), and characterized
by Fiorenza-Karadzhov \cite{fika}, using extrapolation (cf. \cite{kami}) as
follows%
\[
\left\|  f\right\|  _{L_{(q,q^{\prime}}(Q_{n})}\simeq\int_{0}^{1}\left(
\int_{0}^{t}f^{\ast}(s)^{q}ds\right)  ^{1/q}\frac{dt}{t(\log\frac{1}%
{t})^{\frac{1}{2}}}.
\]
The (non homogeneous) result in \cite{FKS} yields%
\begin{equation}
\left\|  f\right\|  _{L_{(q,q^{\prime}}(Q_{n})}\leq C(q)\left(  \left\|
\nabla f\right\|  _{L^{q}(Q_{n})}+\left\|  f\right\|  _{L^{q}(Q_{n})}\right)
,\text{ }f\in C_{0}^{\infty}(Q_{n}). \label{ladesigualdad}%
\end{equation}
Consequently, since%
\[
L_{(q,q^{\prime}}(Q_{n})\subset L^{q}(LogL)^{q/2}(Q_{n}),
\]
this result provides an improvement upon (\ref{launoprima}).

In \cite{mamidim} we showed that the inequality (\ref{ladesigualdad}) is
connected with a different transference principle. We start by reformulating
(\ref{lavieja}) as
\[
\left\Vert \left(  f_{\mu}^{\ast\ast}(t)-f_{\mu}^{\ast}(t)\right)
\chi_{(0,1/2)}(t)G_{\infty}(t)\right\Vert _{\bar{X}}\leq c\left\Vert
G_{\infty}(t)\frac{t}{I(t)}\right\Vert _{L^{\infty}(0,\frac{1}{2})}\left\Vert
\left\vert \nabla f\right\vert _{\mu}^{\ast\ast}\right\Vert _{\bar{X}},
\]
then (\ref{lavieja}) corresponds to the choice $G_{\infty}(t)=(\log\frac{1}%
{t})^{\frac{1}{2}}.$ Now, using the fact that $\frac{I_{\mu}(t)}{t}$ decreases
we see that left hand side of (\ref{latia}) can be minorized as follows,%
\begin{align*}
\left\Vert \left(  f_{\mu}^{\ast\ast}(\cdot)-f_{\mu}^{\ast}(\cdot)\right)
\frac{I_{\mu}(\cdot)}{(\cdot)}\right\Vert _{\bar{X}}  &  \geq\left\Vert
\left(  f_{\mu}^{\ast\ast}(\cdot)-f_{\mu}^{\ast}(\cdot)\right)  \chi
_{(0,t)}(\cdot)\frac{I_{\mu}(\cdot)}{(\cdot)}\right\Vert _{\bar{X}}\\
&  \geq\left\Vert \left(  f_{\mu}^{\ast\ast}(\cdot)-f_{\mu}^{\ast}%
(\cdot)\right)  \chi_{(0,t)}(\cdot)\right\Vert _{\bar{X}}\frac{I_{\mu}(t)}{t}.
\end{align*}
Therefore, we have%
\begin{equation}
\left\Vert \left(  f_{\mu}^{\ast\ast}(\cdot)-f_{\mu}^{\ast}(\cdot)\right)
\chi_{(0,t)}(\cdot)\right\Vert _{\bar{X}}\leq c\frac{t}{I_{\mu}(t)}\left\Vert
\left\vert \nabla f\right\vert _{\mu}^{\ast\ast}\right\Vert _{\bar{X}}.
\label{larusa}%
\end{equation}
Now, if $G_{1}$ is such that $\left(  \int_{0}^{1}G_{1}(t)\frac{t}%
{I(t)}dt\right)  <\infty,$ it follows immediately from (\ref{larusa}) that%
\[
\int_{0}^{1}\left\Vert \left(  f_{\mu}^{\ast\ast}(\cdot)-f_{\mu}^{\ast}%
(\cdot)\right)  \chi_{(0,t)}(\cdot)\right\Vert _{\bar{X}}G(t)dt\leq C\left(
\int_{0}^{1}G_{1}(t)\frac{t}{I(t)}dt\right)  \left\Vert \left\vert \nabla
f\right\vert _{\mu}^{\ast\ast}\right\Vert _{\bar{X}}.
\]
For example, let $G_{1}(t)=\frac{1}{t\left(  \log\frac{1}{t}\right)
^{\frac{1}{2}}},$ and suppose the following (stronger) Gaussian isoperimetric
transference condition is satisfied,
\begin{equation}
\int_{0}^{1}\frac{dt}{I(t)(\log\frac{1}{t})^{\frac{1}{2}}}<\infty.
\label{laconversa}%
\end{equation}
Then, we have (cf. \cite{mamidim})
\begin{align}
\left\Vert (f_{\mu}^{\ast\ast}(\cdot)-f_{\mu}^{\ast}(\cdot))\left(  \log
(\frac{1}{\cdot})\right)  ^{1/2}\right\Vert _{\bar{X}}  &  \leq c\int_{0}%
^{1}\left\Vert \left(  f_{\mu}^{\ast\ast}(s)-f_{\mu}^{\ast}(s)\right)
\chi_{(0,t)}(s)\right\Vert _{\bar{X}}\frac{dt}{t\left(  \log\frac{1}%
{t}\right)  ^{\frac{1}{2}}}\nonumber\\
&  \leq c\left(  \int_{0}^{1}\frac{dt}{I(t)(\log\frac{1}{t})^{\frac{1}{2}}%
}\right)  \left\Vert \left\vert \nabla f\right\vert ^{\ast}\right\Vert
_{\bar{X}}.\nonumber
\end{align}

Let us show a concrete application. Let
\[
I_{n}(t)=n\left(  \gamma_{n}\right)  ^{1/n}t^{1-1/n},
\]
where $\gamma_{n}=\frac{\pi^{n/2}}{\Gamma(1+n/2)}$ is the measure of the unit
ball in $\mathbb{R}^{n}$ (i.e. $I_{n}(t)$ is the isoperimetric profile
associated to $\mathbb{R}^{n}),$ and consider the function%
\[
G_{1}(t)=\frac{1}{t\sqrt{\ln\left(  \frac{1}{t}\right)  }},t\in(0,1).
\]
Then,%
\begin{align*}
\int_{0}^{1}\frac{t}{tI_{n}(t)}G_{1}(t)dt  &  =\frac{1}{n\left(  \gamma
_{n}\right)  ^{1/n}}\int_{0}^{1}t^{1/n}\frac{dt}{t\left(  \ln\frac{1}%
{t}\right)  ^{\frac{1}{2}}}\\
&  =\frac{1}{n\left(  \gamma_{n}\right)  ^{1/n}}\int_{0}^{\infty}z^{-\frac
{1}{2}}e^{-z/n}dz\\
&  =\frac{\sqrt{\pi}n^{\frac{1}{2}}}{n\left(  \gamma_{n}\right)  ^{1/n}}\\
&  =\frac{\Gamma(1+\frac{n}{2})^{1/n}}{n^{\frac{1}{2}}}.\\
&  =\left(  \frac{n}{2}\right)  ^{1/n}\frac{\Gamma(\frac{n}{2})^{1/n}%
}{n^{\frac{1}{2}}}\\
&  \leq\frac{1}{\sqrt{2}}\left(  \frac{n}{2}\right)  ^{1/n}\\
&  \leq c.
\end{align*}

Thus,
\[
\sup_{n}\int_{0}^{1}\frac{dt}{I_{n}(t)(\log\frac{1}{t})^{\frac{1}{2}}}%
<\infty.
\]
As a consequence the following dimensionless Sobolev inequality holds,
\[
\int_{0}^{1}\left\|  \left(  f^{\ast\ast}(\cdot)-f^{\ast}(t)\right)
\chi_{\lbrack0,t)}(\cdot)\right\|  _{\bar{X}}\frac{dt}{t(\log\frac{1}%
{t})^{1/2}}\leq C\left\|  \left|  \nabla f\right|  ^{\ast\ast}\right\|
_{\bar{X}}.
\]
For $X=L^{q}$ and for $\Omega=Q_{n}$ this gives the result of \cite{FKS}.

\section{Rearrangement inequalities of Garsia-Rodemich type and Morrey's
theorem}

To complement the results of previous sections we now consider fractional
inequalities and the Morrey-Sobolev embedding theorem. Here we follow
\cite{mamiaster} where the reader will find a complete treatment together with
many applications.

Let us describe model results that influenced our development in these directions.

\begin{example}
\label{ejemplo1}For all $f\in X(\mathbb{R}^{n})+\dot{W}_{X}^{1}(
\mathbb{R}^{n}),$ we have (cf. \cite{mamipams})%
\begin{equation}
f^{\ast\ast}(t)-f^{\ast}(t)\leq c_{n}\frac{\omega_{X}\left(  t^{1/n},f\right)
}{\phi_{X}(t)},\text{ }t>0, \label{tres}%
\end{equation}
where $X=X(\mathbb{R}^{n})$ is a rearrangement invariant space on
$\mathbb{R}^{n},$ $\phi_{X}(t)=\left\|  \chi_{A}\right\|  _{X},$ with $\left|
A\right|  =t,$ is the fundamental function of $X,$ and $\omega_{X}$ be the
modulus of continuity associated with $X$:
\[
\omega_{X}\left(  t,g\right)  =\sup_{\left|  h\right|  \leq t}\left\|
g(\cdot+h)-g(\cdot)\right\|  _{X},\text{ for }g\in X.
\]
The inequality (\ref{tres}) can be formulated and proved on different levels
of generality on the spaces, the domains; and indeed have a long history: We
refer to \cite{kolyada}, \cite{jonen}, \cite{BS}, \cite{mamipams} and the
references therein.
\end{example}

The associated questions to Example \ref{ejemplo1} are: What is the
corresponding Gaussian result? More generally: What is the metric version?
What is the role of dimension? What is the connection with isoperimetry? Even
to formulate metric results we need to develop suitable tools. For example:
What is an appropriate replacement for the modulus of continuity?

We shall need the following definition: Consider a connected, measure metric
spaces $\left(  \Omega,d,\mu\right)  $ equipped with a finite Borel measure
$\mu$. For measurable functions $u:\Omega\rightarrow\mathbb{R}$, the signed
\textbf{decreasing rearrangement} of $u,$ which we denote by $u_{\mu}^{s},$ is
the right-continuous non-increasing function from $[0,\mu(\Omega))$ into
$\mathbb{R}$ that is equimeasurable with $u;$ i.e. $u_{\mu}^{s}$ satisfies
\[
\mu\{x\in{\Omega}:u(x)>t\}=m(\left\{  z\in\lbrack0,\mu(\Omega)):u_{\mu}%
^{s}(z)>t\right\}  )\text{ , \ }t\in\mathbb{R}%
\]
The maximal average$\ $of $u_{\mu}^{ss}$ is defined by
\[
u_{\mu}^{ss}(t)=\frac{1}{t}\int_{0}^{t}u_{\mu}^{s}(z)dz,\text{ }(t>0).
\]

\begin{example}
\label{ejemplo2}In closely related work Garsia and his collaborators (cf.
\cite{garro}, \cite{garsiaind}, \cite{garsia} and the references therein)
obtained related inequalities. For example, in \cite{garsiaind} and
\cite{garsia} for functions defined on the unit $n-$cube, and using signed
rearrangements with respect the Lebesgue measure
\begin{equation}
\left.
\begin{array}
[c]{c}%
f^{s}(x)-f^{s}(1/2)\\
f^{s}(1/2)-f^{s}(1-x)
\end{array}
\right\}  \leq c\int_{x}^{1}\frac{\omega_{L^{p}}(t^{1/n},f)}{t^{1/p}}\frac
{dt}{t},x\in(0,\frac{1}{2}], \label{degarsia}%
\end{equation}
where $\omega_{L^{p}}(t,f)$ the $L^{p}$ modulus of continuity. The extension
from dimension $1$ to dimension $n$ in these works was done through highly non
trivial combinatorial inequalities (cf. \cite{garsia}).
\end{example}

In this context we can ask similar questions to those posed in Example
\ref{ejemplo1}.

\begin{example}
\label{ejemplo3}In the work of Garsia and his collaborators one also finds
another interesting use of rearrangement inequalities to extract continuity.
Through a change of scale argument, inequalities on rearrangements were used,
for example, to prove versions of Morrey's Sobolev theorem. For example, in
the one dimensional case (cf. \cite{garsiaind}, \cite{garro}, \cite{garsia})
we have%
\[
\left\vert f(x)-f(y)\right\vert \leq2c\int_{0}^{\left\vert x-y\right\vert
}\frac{\omega_{L^{p}}(t,f)}{t^{1/p}}\frac{dt}{t};\text{ }x,y\in\lbrack0,1].
\]
We now adapt the change of scale argument of Garsia et. al. combined with the
rearrangement inequality (\ref{metricas}) in the context of the unit cube to
prove a version of Morrey's theorem. First, let us stipulate (cf.
\cite{mamiaster}) that one can rewrite the (\ref{metricas}) in terms of signed
rearrangements (i.e. we rearrange $f$ rather than its absolute value). Suppose
that $p>n,$ and let $f\in W_{1}^{1}(0,1)^{n}.$ Starting with the fundamental
theorem of calculus%
\[
f^{ss}(0)-f^{ss}(1)=\int_{0}^{1}\left(  f^{ss}(t)-f^{s}(t)\right)  \frac
{dt}{t}%
\]
and the corresponding version of (\ref{metricas})
\[
f^{ss}(t)-f^{s}(t)\leq c_{n}\frac{t}{\min(t,1-t)^{1-1/n}}\left\vert \nabla
f\right\vert ^{\ast\ast}(t),\text{ \ }0<t<1,
\]
we see that%
\begin{align}
f^{ss}(0)-f^{ss}(1) &  \leq c_{n}\int_{0}^{1}\left\vert \nabla f\right\vert
^{\ast\ast}(t)\frac{dt}{\min(t,1-t)^{1-1/n}}\nonumber\\
&  \leq c_{n,p}\left\Vert \left\vert \nabla f\right\vert \right\Vert _{L^{p}%
}\left\Vert \frac{1}{\min(t,1-t)^{1-1/n}}\right\Vert _{L^{p^{\prime}}%
(0,1)}\text{ (by H\"{o}lder's inequality)}\nonumber\\
&  =C_{n.p}\left\Vert \left\vert \nabla f\right\vert \right\Vert _{L^{p}%
}.\nonumber
\end{align}
This works of course because for $p>n,\left\Vert \frac{1}{\min(t,1-t)^{1-1/n}%
}\right\Vert _{L^{p^{\prime}}}<\infty.$ Now (recall we are working with signed
rearrangements) since $f^{ss}(0)=ess\sup_{x\in(0,1)^{n}}f,$ and $f^{ss}%
(1)=\int_{0}^{1}f,$ it thus follows that
\begin{equation}
ess\sup_{x\in(0,1)^{n}}f-\int_{0}^{1}f=f^{ss}(0)-f^{ss}(1)\leq C_{n.p}%
\left\Vert \left\vert \nabla f\right\vert \right\Vert _{L^{p}}.\label{vale}%
\end{equation}
Applying (\ref{vale}) now to $-f$ \ yields
\begin{equation}
\int_{0}^{1}f-ess\inf_{x\in(0,1)^{n}}f\leq C_{n.p}\left\Vert \left\vert \nabla
f\right\vert \right\Vert _{L^{p}}.\label{vale1}%
\end{equation}
Therefore, adding (\ref{vale}) and (\ref{vale1}) we obtain
\[
Osc(f;(0,1)^{n}):=ess\sup_{x\in(0,1)^{n}}f-ess\inf_{x\in(0,1)^{n}}%
f\leq2C_{n.p}\left\Vert \left\vert \nabla f\right\vert \right\Vert _{L^{p}}.
\]
Now, we scale: Apply the previous argument to the restriction of $f$ to a
subcube $Q,$ $f\chi_{Q}$. We obtain,%
\[
Osc(f;Q)\leq c_{n,p}\left\Vert \frac{t}{\min(t,\left\vert Q\right\vert
-t)^{1-1/n}}\right\Vert _{L^{p^{\prime}}(0,|Q|)}\left\Vert \left\vert \nabla
f\right\vert \right\Vert _{L^{p}(Q)}.
\]
By computation, it is now easy to see that we will have
\[
\left\vert f(y)-f(z)\right\vert \leq c_{n,p}\left\vert y-z\right\vert
^{(1-\frac{n}{p})}\left\Vert \left\vert \nabla f\right\vert \right\Vert
_{p},\text{a.e. }y,z.
\]

\end{example}

To carry out our program in metric measure spaces $(\Omega,d,\mu)$ we first
need to find a replacement for the modulus of continuity. In \cite{mamiaster}
we use the Peetre\footnote{Independently, and in parallel, A. Calder\'{o}n and
his student Oklander defined and studied the $K-$functional, and real
interpolation, e.g. in Oklander's thesis at the University of Chicago (cf.
\cite{okl}).} $K-$functional defined by%
\begin{align*}
K(t,f;X(\Omega),S_{X}(\Omega))  &  :=\\
\inf\{\left\Vert f-g\right\Vert _{X(\Omega)}+t\left\Vert \left\vert \nabla
g\right\vert \right\Vert _{X(\Omega)}  &  :g\in S_{X}(\Omega)\},
\end{align*}
where for a rearrangement invariant space $X(\Omega)$, $S_{X}(\Omega)=\{f\in
Lip(\Omega):$ $\left\Vert \left\vert \nabla f\right\vert \right\Vert
_{X(\Omega)}<\infty\}.$ In the classical setting we have (cf. \cite{BS},
\cite{bl})
\begin{align*}
K(t,f;X(\mathbb{R}^{n}),\dot{W}_{X}^{1}(\mathbb{R}^{n}))  &  :=\inf
\{\left\Vert f-g\right\Vert _{X}+t\left\Vert \left\vert \nabla g\right\vert
\right\Vert _{X}:g\in\dot{W}_{X}^{1}(\mathbb{R}^{n})\}\\
&  \simeq\omega_{X}(t,f).
\end{align*}
We can thus think of $K(t,f;X(\Omega),S_{X}(\Omega))$ as generalized
\textquotedblleft modulus of continuity\textquotedblright, and we have (cf.
\cite{mamiaster} and also \cite{mas})

\begin{theorem}
(cf. \cite{mamiaster}) Let $\left(  \Omega,d,\mu\right)  $ be a metric measure
space that satisfies our standard assumptions. Then,

(i) For all rearrangement invariant spaces $X(\Omega)$, and for all $f\in
X(\Omega)+S_{X}(\Omega),$
\begin{equation}
f_{\mu}^{ss}(t)-f_{\mu}^{s}(t)\leq16\frac{K\left(  \frac{t}{I_{\Omega}%
(t)},f;X(\Omega),S_{X}(\Omega)\right)  }{\phi_{X}(t)},\text{ }t\in
(0,\mu(\Omega)/2), \label{desiK}%
\end{equation}

\begin{equation}
\left(  f-f_{\Omega}\right)  _{\mu}^{\ast\ast}(t)-\left(  f-f_{\Omega}\right)
_{\mu}^{\ast}(t)\leq16\frac{K\left(  \frac{t}{I_{\Omega}(t)},f;X(\Omega
),S_{X}(\Omega)\right)  }{\phi_{X}(t)},\text{ }t\in(0,\mu(\Omega)),
\label{desiK9}%
\end{equation}
where
\begin{equation}
f_{\Omega}=\frac{1}{\mu(\Omega)}\int_{\Omega}fd\mu. \label{average}%
\end{equation}

(ii) Conversely, suppose that $G:(0,\mu(\Omega))\rightarrow\mathbb{R}_{+},$ is
a continuous function, which is concave and symmetric around $\mu(\Omega)/2,$
and that there exists a constant $c>0$ such that\footnote{In other words we
assume that (\ref{desiK}) holds for $X=L^{1}(\Omega),$ and with $\frac
{t}{G(t)}$ replacing $\frac{t}{I_{\Omega}(t)}.$} for all $f\in X(\Omega
)+S_{X}(\Omega),$%
\[
f_{\mu}^{ss}(t)-f_{\mu}^{s}(t)\leq c\frac{K\left(  \frac{t}{G(t)}%
,f;X(\Omega),S_{X}(\Omega)\right)  }{t},\text{ }t\in(0,\mu(\Omega)/2).
\]
Then, there exists a constant $c>0$ such that for all $t\in(0,\mu(\Omega)),$%
\[
G(t)\leq cI_{\Omega}(t).
\]

\end{theorem}

Following the analysis of \cite[Chapter 10]{mamiaster}, one can show that for
a metric probability space $\left(  \Omega,d,\mu\right)  ,$ that satisfies our
standard assumptions, the Garsia inequality (\ref{degarsia}) takes the
following form,
\begin{equation}
\left.
\begin{array}
[c]{c}%
f^{s}(x)-f^{s}(1/2)\\
f^{s}(1/2)-f^{s}(1-x)
\end{array}
\right\}  \leq c\int_{x}^{1}\frac{K\left(  \frac{t}{I_{\Omega}(t)}%
,f;X(\Omega),S_{X}(\Omega)\right)  }{\phi_{X}(t)}\frac{dt}{t},x\in(0,\frac
{1}{2}].\label{degarsia2}%
\end{equation}

For the scaling argument we outlined above we need an extra assumption. We say
that an \textbf{isoperimetric inequality relative to} $G$ holds, if there
exists a positive constant $C_{G}$ such that
\[
I_{G}(s)\geq C_{G}\min(I_{\Omega}(s),I_{\Omega}(\mu(G)-s)).
\]
We say that the metric measure space $\left(  \Omega,d,\mu\right)  $ has the
\textbf{uniform} \textbf{relative isoperimetric property,} if there exist
$C>0,$ $\delta>0,$ such that for any $x\in\Omega,$ and any open ball
$B_{\alpha}(x)$ centered on $x,$ with $\mu(B_{\alpha}(x))=\alpha$
($0<\alpha<\delta),$ the \textbf{relative isoperimetric profile }%
$I_{B_{\alpha}(x)}$ satisfies:
\[
I_{B_{\alpha}(x)}(s)\geq C\min(I_{\Omega}(s),I_{\Omega}(\alpha-s)),\text{
\ \ }0<s<\alpha.
\]

Then we have (cf. \cite[Chapter 4]{mamiaster})

\begin{theorem}
Let $\left(  \Omega,d,\mu\right)  $ be a metric measure space that satisfies
our standard assumptions and, moreover, has the relative uniform isoperimetric
property\textbf{. }Let $X$ be a r.i. space in $\Omega$ such that
\[
\left\Vert \frac{1}{I_{\Omega}(s)}\right\Vert _{\bar{X}^{^{\prime}}}<\infty.
\]
Then, if $f\in X+S_{X}(\Omega)$ satisfies
\[
\int_{0}^{\mu(\Omega)}\frac{K\left(  \phi_{X}(t)\left\Vert \frac{1}{I_{\Omega
}(s)}\chi_{(0,t)}(s)\right\Vert _{\bar{X}^{^{\prime}}},f;X,S_{X}%
(\Omega)\right)  }{\phi_{_{X}}(t)}\frac{dt}{t}<\infty,
\]
it follows that $f$ is essentially bounded and essentially continuous.
\end{theorem}

For applications we need to have explicit computations of the relevant
$K-$functional. We refer to \cite{ditzian}, \cite{ditiv}, \cite{dil} for a
treatment of $K-$functionals in one dimensional Gaussian measure. For other
relevant computations of $K-$functionals see the appendix of \cite{mamiaster}.

For further developments we must refer to \cite{mamiaster}. In connection with
this section we should also mention \cite{xiao} and the references therein.

\section{Higher Order Inequalities}

In this section we work with Euclidean domains $\Omega$ and measure spaces of
the form $d\mu(x)=w(x)dx,$ and we formulate higher order pointwise
inequalities by iteration\footnote{The iteration of Sobolev inequalities is
not a new idea (cf. \cite[Corollary 6.9/1 in page 379 and Theorem 7.6.5 in
page 430.]{maz}), the novelty here, if any, is the iteration of pointwise
rearrangement inequalities} that extend our previous work in \cite{mamijmaa}
and \cite{mp}. The basic inequality here reads as follows. We assume that
$\int d\mu(x)<\infty,$ and we let $I_{\mu}:=I,$ then for sufficiently smooth
$f$\ we have that for $k\geq2$, (cf. \cite{mamiitera}))%

\begin{align*}
f^{\ast\ast}(t)-f^{\ast}(t)  &  \leq\frac{1}{(k-1)!}\frac{t}{I(t)}\int
_{t}^{1/2}\left\vert d^{k}f\right\vert ^{\ast\ast}(u)\left(  \frac{1}%
{I(u)}\left(  \int_{t}^{u}\frac{dz}{I(z)}\right)  ^{k-1}\right)  du\\
&  +\frac{1}{(k-1)!}\sum_{j=1}^{k-1}\left(  \frac{t}{I(t)}\left(  \int
_{t}^{1/2}\frac{dz}{I(z)}\right)  ^{k-j-1}\right)  \left\Vert \left\vert
d^{k-j}f\right\vert \right\Vert _{1},\text{ }0<t<1/2.
\end{align*}
As a corollary we obtain,%

\[
f^{\ast\ast}(t)\leq c_{k}\int_{t}^{1/2}\left\vert d^{k}f\right\vert ^{\ast
\ast}(s)\frac{1}{I(s)}\left(  \int_{t}^{s}\frac{dz}{I(z)}\right)  ^{k}%
ds+\sum_{j=1}^{k}c_{j}(t)\left\Vert \left\vert d^{k-j}f\right\vert \right\Vert
_{1},\text{ }0<t<1/2;
\]
where $c_{j}(t)=\left(  \int_{t}^{1/2}\frac{dz}{I(z)}\right)  ^{k-j}.$

In particular, one can easily show that this result implies recent iterated
inequalities that appeared in \cite{slavikova}.

\begin{theorem}
Let $X,Y$ r.i spaces. Suppose that the operator $T$ defined by%
\[
Tf(t)=\int_{t}^{1/2}\frac{f(u)}{I(u)}\left(  \int_{t}^{u}\frac{dz}%
{I(z)}\right)  ^{k-1}du,
\]
is bounded from $X$ to $Y.$ Then%
\[
\left\Vert u\right\Vert _{Y}\leq c\left(  \left\Vert \left\vert d^{k}%
f\right\vert ^{\ast\ast}\right\Vert _{X}+\sum_{j=1}^{k}\left\Vert \left\vert
d^{k-j}f\right\vert \right\Vert _{1}\right)  .
\]

\end{theorem}

\begin{acknowledgement}
We are grateful to E. Milman for a number of useful comments that helped
improve the presentation.
\end{acknowledgement}


\begin{thebibliography}{999}                                                                                              %
\bibitem {adams}D. R. Adams and L. I. Hedberg, \textsl{Function spaces and
potential theory}, Grund. Math. \textbf{314}, Springer-Verlag, 1999.

\bibitem {allen}G. D. Allen, \textsl{Locally Continuous Operators II}, Indiana
Univ. Math. J. \textbf{38} (1989), 711--743.

\bibitem {alm}F. Almgren and E. Lieb, \textsl{Symmetric Decreasing
Rearrangement is sometimes continuous,} J. Amer. Math. Soc. \textbf{2} (1989), 683-773.

\bibitem {alv}A. Alvino, G. Trombetti and P. L. Lions, \textsl{On optimization
problems with prescribed rearrangements}, Nonlinear Anal. \textbf{13} (1989), 185-220.

\bibitem {bakr}D. Bakry, T. Coulhon, M. Ledoux and L. Saloff-Coste,
\textsl{Sobolev inequalities in disguise}, Indiana Univ. Math. J. \textbf{44}
(1995) 1033--1074.

\bibitem {bart1}F. Barthe, \textsl{Levels of concentration between exponential
and Gaussian}, Ann. Fac. Sci. Toulouse Math. \textbf{10} (2001), 393--404.

\bibitem {bart}F. Barthe,\textsl{ Log-concave and spherical models in
isoperimetry}, Geom. Funct. Anal. \textbf{12} (2002), 32--55.

\bibitem {BCR}F. Barthe, P. Cattiaux and C. Roberto, \textsl{Isoperimetry
between exponential and Gaussian}, \textsl{Orlicz hyper-contractivity and
isoperimetry}, Rev. Mat. Iber. \textbf{22} (2006), 993--1067.

\bibitem {BCR1}F. Barthe, P. Cattiaux and C. Roberto, \textsl{Isoperimetry
between exponential and Gaussian}, Electronic Journal of Probability
\textbf{12} (2007), 1212--1237.

\bibitem {bmr}J. Bastero, M. Milman and F. Ruiz,\textsl{ On the connection
between weighted norm inequalities, commutators and real interpolation},
preprint, Sem A. Galdeano, 1996.

\bibitem {bmrind}J. Bastero, M. Milman and F. Ruiz, \textsl{A note on
L(}$\infty$\textsl{, q) spaces and Sobolev embeddings}, Indiana Mathematics
Journal \textbf{52} (2003), 1215-1230.

\bibitem {bayle}V. Bayle, \textsl{Propri\'{e}t\'{e}s de concavit\'{e} du
profil isop\'{e}rim\'{e}trique et applications}. Ph.D. thesis, Institut Joseph
Fourier, Grenoble (2004).

\bibitem {be2}W. Beckner and M. Persson, \textsl{On sharp Sobolev embedding
and the logarithmic Sobolev inequality}, Bull. London Math. Soc. \textbf{30}
(1998), 80-84

\bibitem {bds}C. Bennett, R. DeVore, and R. Sharpley, \textsl{Weak-L}%
$^{\infty}$\textsl{ and BMO}, Annals of Math. \textbf{113} (1981), 601-611.

\bibitem {BS}C. Bennett and R. Sharpley, \textsl{Interpolation of Operators},
Academic Press, Boston, 1988.

\bibitem {bl}J. Bergh and J. L\"{o}fstr\"{o}m, \textsl{Interpolation spaces.
An introduction}, Springer\-Verlag, Berlin-Heidelberg-New York,
1976{\footnotesize .}

\bibitem {Bob1}S. G. Bobkov. \textsl{Isoperimetric and analytic inequalities
for log-concave probability measures}, Ann. Probab. \textbf{27} (1999), 1903--1921.

\bibitem {bohou}S.G. Bobkov and C. Houdr\'{e}, \textsl{Some connections
between isoperimetric and Sobolev-type inequalities}, Mem. Amer. Math. Soc.
\textbf{129} (1997).

\bibitem {bobzeg}S. G. Bobkov and B. Zegarlinski, \textsl{Entropy bounds and
isoperimetry}, Mem. Amer. Math. Soc. \textbf{176} (2005).

\bibitem {bobz}S. G. Bobkov and B. Zegarlinski,\textit{ }\textsl{Distributions
with slow tails and ergodicity of Markov semigroups in infinite dimensions},
in A. Laptev (ed), Around the research of Vladimir Maz'ya I: Function Spaces,
Springer, 2010, pp 13-79.

\bibitem {Bor1}C. Borell, \textsl{The Ehrhard inequality}, C. R. Math. Acad.
Sci. Paris \textbf{337} (2003), 663--666.

\bibitem {bo}C. Borell, \textsl{The Brunn-Minkowski inequality in Gauss
space}, Invent. Math. \textbf{30} (1975), 207-216.

\bibitem {bw}H. Brezis and S. Wainger, \textsl{A note on limiting cases of
Sobolev embeddings and convolution inequalities}, Comm. Partial Diff. Eq.
\textbf{5} (1980), 773-789.

\bibitem {cal}A. P. Calder\'{o}n, \textsl{Spaces between }$L^{1}$\textsl{ and
}$L^{\infty}$\textsl{ and the theorem of Marcinkiewicz}, Studia Math.
\textbf{26} (1966), 273-299.

\bibitem {Ciapick}A. Cianchi and L. Pick, \textsl{Optimal Gaussian Sobolev
embeddings}, J. Funct. Anal. \textbf{256} (2009), 3588--3642.

\bibitem {cou}T. Coulhon, \textsl{Espaces de Lipschitz et in\'{e}galit\'{e}s
de Poincar\'{e}}, J. Funct. Anal. \textbf{136} (1996), 81-113.

\bibitem {cou2}T. Coulhon, \textsl{Dimensions at infinity for Riemannian
manifolds}, Pot. Anal. \textbf{4} (1995), 335-344.

\bibitem {cou1}T. Coulhon, \textsl{Heat kernel and isoperimetry on non-compact
Riemmanian manifolds}, Contemporary Mathematics \textbf{338} (2003), 65-99.

\bibitem {cwijami}M. Cwikel, B. Jawerth and M. Milman, \textsl{A note on
extrapolation of inequalities}, preprint, 2010.

\bibitem {cwp}M. Cwikel and E. Pustylnik, \textsl{Sobolev type embeddings in
the limiting case}, J. Fourier Anal. Appl. \textbf{4} (1998), 433--446.

\bibitem {ditzian}Z. Ditzian and V. Totik, \textsl{Moduli of Smoothness},
Springer-Verlag, New York, 1987.

\bibitem {ditiv}Z. Ditzian and K. G. Ivanov, \textsl{Strong converse
inequalities}, J. D'Analise Math. \textbf{61} (1993), 61-111.

\bibitem {dil}Z. Ditzian and D. S. Lubinsky, \textsl{Jackson and Smoothness
theorems for Freud weights in }$L_{p}$ $(0<p<\infty),$ Constr. Approx.
\textbf{13} (1997), 99-152.

\bibitem {edmund}D. E. Edmunds and W. D. Evans, \textsl{Hardy operators,
function spaces and embeddings}, Springer-Verlag, Berlin, 2004.

\bibitem {er}A. Ehrhard, \textsl{Sym\'{e}trisation dans l'espace de Gauss},
Math. Scand. \textbf{53} (1983), 281-301.

\bibitem {er1}A. Ehrhard, \textsl{Sur l'in\'{e}galit\'{e} de Sobolev
logarithmique de Gross}, S\'{e}minaire de Probabilit\'{e}s XVI I I, Lecture
Notes in Math. 1059, 194-196, Springer-Verlag, 1984.

\bibitem {er2}A. Ehrhard, \textsl{In\'{e}galit\'{e}s isop\'{e}rim\'{e}triques
et int\'{e}grales de Dirichlet gaussiennes}, Ann. Scient. Ec. Norm. Sup.
\textbf{17} (1984), 317-332.

\bibitem {fio1}A. Fiorenza, \textsl{Duality and reflexivity in grand Lebesgue
spaces}, Collect. Math. \textbf{51}, (2000), 131-148.

\bibitem {fika}A. Fiorenza and G. E. Karadzhov, \textsl{Grand and Small
Lebesgue Spaces and their analogs}, Z. Anal. Anwendungen, \textbf{23} (2004), 657-681.

\bibitem {FKS}A. Fiorenza, M. Krbec and H. J. Schmeisser, \textsl{An
improvement of dimension-free Sobolev imbeddings in r.i. spaces}, preprint.

\bibitem {carlo}L. Fontana and C. Morpurgo, \textsl{Optimal limiting
embeddings for }$\Delta$\textsl{-reduced Sobolev spaces in }$L^{1},$ Ann. de
l'Inst. Henri Poincar\'{e} (C) Non Linear Analysis, to appear.

\bibitem {gallot}S. Gallot, \textsl{In\'{e}galit\'{e}s
isop\'{e}rim\'{e}triques et analytiques sur les vari\'{e}t\'{e}s
Riemanniennes}, Ast\'{e}risque No. \textbf{163-164} (1988), 31--91.

\bibitem {garsia}A. M. Garsia, \textsl{Combinatorial inequalities and
smoothness of functions}, Bull. Amer. Math. Soc. \textbf{82} (1976), 157-170.

\bibitem {garsiaind}A. M. Garsia, \textsl{A remarkable inequality and the
uniform convergence of Fourier series}, Indiana Univ. Math. J. \textbf{25}
(1976), 85-102.

\bibitem {garro}A. Garsia and E. Rodemich, \textsl{Monotonicity of certain
functionals under rearrangement}, Ann. Inst. Fourier \textbf{24} (1974), 67-116.

\bibitem {griebel}M. Griebel, \textsl{Sparse grids and related approximation
schemes for higher dimensional problems}, in L. Pardo, A. Pinkus, E. Suli, and
M. Todd, editors, Foundations of Computational Mathematics (FoCM05),
Santander, pp 106-161, Cambridge University Press, 2006.

\bibitem {gro}L. Gross, \textsl{Logarithmic Sobolev inequalities}, Amer. J.
Math. \textbf{97} (1975), 1061--1083.

\bibitem {haj}P. Hajlasz,\textsl{\ Sobolev inequalities, truncation method,
and John domains}, Papers in Analysis, Rep. Univ. Jyv\"{a}skyl\"{a} Dep. Math.
Stat. \textbf{83}, Univ. Jyv\"{a}skyl\"{a}, Jyv\"{a}skyl\"{a}, 2001, pp 109--126.

\bibitem {HK2}P. Hajlasz and P. Koskela,\textsl{ Sobolev met Poincar\'{e}},
Mem. Amer. Math. Soc. \textbf{145} (2000), 101 pages.

\bibitem {hansson}K. Hansson, \textsl{Imbedding theorems of Sobolev type in
potential theory}, Math Scand \textbf{45} (1979), 77-102.

\bibitem {hedb}L. I. Hedberg, \textsl{On Maz'ya's work in potential theory and
the theory of function spaces}, The Maz'ya anniversary collection, Vol. 1
(Rostock, 1998), 7--16, Oper. Theory Adv. Appl. \textbf{109}, Birkh\"{a}user,
Basel, 1999.

\bibitem {he}J. Heinonen, \textsl{Lectures on Analysis on metric spaces,
}Lecture Notes, 1996, University of Michigan\textsl{.}

\bibitem {holeem}C. Houdre, M. Ledoux, E. Milman and M. Milman,
\textsl{Concentration, functional inequalities and isoperimetry}, Contemporary
Mathematics \textbf{545}, 2011.

\bibitem {iw}T. Iwaniec and C. Sbordone, \textsl{On the integrability of the
Jacobian under minimal hypotheses}, Arch. Rational Mech. Anal. \textbf{119}
(1992), 129-143.

\bibitem {jm}B. Jawerth and M. Milman, \textsl{Extrapolation theory with
applications}, Mem. Amer. Math. Soc. \textbf{89} (1991), no. 440.

\bibitem {jm2}B. Jawerth and M. Milman,\textsl{ Interpolation of Weak Type
Spaces}, Math. Z. \textbf{201} (1989), 509 - 520)

\bibitem {jonen}H. Johnen and K. Scherer, \textsl{On the equivalence of the
K-functional and moduli of continuity and some applications}, in Constructive
theory of functions of several variables, Lecture Notes in Math. \textbf{571},
pp. 119-140, Springer, Berlin, 1977.

\bibitem {kara}G. E. Karadzhov and Q. Mehmood, \textsl{Optimal Regularity
Properties of the Generalized Sobolev Spaces}, J. Funct. Spaces Appl. (2013).

\bibitem {kami}G. E. Karadzhov and M. Milman, \textsl{Extrapolation Theory:
New Results and Applications, }J. Approx. Th. \textbf{133} (2005), 38-99.

\bibitem {kesavan}S. Kesavan, \textsl{Symmetrization and applications}, World
Scientific, 2006.

\bibitem {kolyada}V.I. Kolyada, \textsl{Rearrangements of functions and
embedding theorems}, Uspekhi Mat. Nauk \textbf{44} (1989) 61-95; transl. in:
Russian Math. Surveys \textbf{44} (1989), 73-117.

\bibitem {krb1}M. Krbec and H. J. Schmeisser, \textsl{On dimension-free
Sobolev imbeddings I}, J. Math. Anal. Appl. \textbf{387} (2012), 114-125.

\bibitem {krb2}M. Krbec and H. J. Schmeisser, \textsl{On dimension-free
Sobolev imbeddings II}, Rev. Mat. Complutense \textbf{25} (2012), 247-265.

\bibitem {le1}M. Ledoux,\textsl{ Isoperimetry and Gaussian Analysis}, Ecole
d'Et\'{e} de Probabilit\'{e}s de Saint-Flour 1994, Springer Lecture Notes
1648, pp 165-294, Springer-Verlag, 1996.

\bibitem {led}M. Ledoux, \textsl{Isop\'{e}rim\'{e}trie et in\'{e}galit\'{e}es
de Sobolev logarithmiques gaussiennes}, C. R. Acad. Sci. Paris Ser. I Math.
\textbf{306} (1988), 79-92.

\bibitem {ledouxbk}M. Ledoux, \textsl{The Concentration of Measure
Phenomenon}, Math. Surveys \textbf{89}, Amer. Math. Soc., 2001.

\bibitem {ledcon}M. Ledoux, \textsl{From concentration to isoperimetry:
Semigroup proofs}, Cont. Math. \textbf{545} (2011), 155-166.

\bibitem {leoni}G. Leoni, \textsl{A first course in Sobolev spaces}, Grad.
Studies in Math. \textbf{105}, Amer. Math. Soc. 2009.

\bibitem {mapi}J. Mal\'{y} and L. Pick, \textsl{An elementary proof of Sharp
Sobolev embeddings}, Proc. Amer. Math. Soc. \textbf{130} (2002), 555-563.

\bibitem {mamipams}J. Mart\'{\i}n and M. Milman, \textsl{Symmetrization
inequalities and Sobolev embeddings}, Proc. Amer. Math. Soc. \textbf{134}
(2006), 2335-2347.

\bibitem {mamijmaa}J. Mart\'{\i}n and M. Milman, \textsl{Higher-order
symmetrization inequalities and applications}, J. Math. Anal. Appl.
\textbf{330} (2007), 91-113.

\bibitem {mami07}J. Martin and M. Milman, \textsl{A note on Sobolev
inequalities and limits of Lorentz spaces}, Contemp. Math. \textbf{445}
(2007), 237-245.

\bibitem {mamijfa}J. Martin and M. Milman, \textsl{Isoperimetry and
Symmetrization for Logarithmic Sobolev inequalities}, Journal of Functional
Analysis \textbf{256} (2009), 149-178.

\bibitem {mamiarxiv}J. Martin and M. Milman, \textsl{Addendum to Isoperimetry
and Symmetrization for Logarithmic Sobolev inequalities,} (arXiv:0901.1839)

\bibitem {mamicomptes}J. Martin and M. Milman, \textsl{Isoperimetry and
Symmetrization for Sobolev spaces on metric spaces}, Comptes Rendus Math.
\textbf{347} (2009), 627--630.

\bibitem {mamimazya}J. Mart\'{\i}n and M. Milman, \textsl{Isoperimetric Hardy
type and Poincar\'{e} inequalities on metric spaces}, In: Around the Research
of Vladimir Maz'ya I. Function Spaces - Ari Laptev (Ed.) International
Mathematical Series, Springer 11 (2010), 285--298.

\bibitem {mamiadv}J. Martin and M. Milman, \textsl{Pointwise Symmetrization
Inequalities for Sobolev functions and applications}, Adv. Math. \textbf{225}
(2010), 121-199.

\bibitem {mamicon}J. Martin and M. Milman, \textsl{Sobolev inequalities,
rearrangements, isoperimetry and interpolation spaces}, Contemp. Math.
\textbf{545} (2011), 167-193.

\bibitem {mamiaster}J. Martin and M. Milman, \textsl{Fractional Sobolev
inequalities: symmetrization, isoperimetry and interpolation}, submitted (arXiv:1205.1231).

\bibitem {mamicoulhon}J. Martin and M. Milman, \textsl{A note on Coulhon type
inequalities}, to appear in Proc. Amer. Math. Soc. (arXiv:1206.1584)

\bibitem {mamidim}J. Martin and M. Milman, \textsl{Integral isoperimetric
transference and dimensionless Sobolev inequalities}, submitted (arXiv:1309.1980).

\bibitem {mamiitera}J. Martin and M. Milman, \textsl{A note on iterated
Sobolev inequalities involving the isoperimetric profile}, preprint.

\bibitem {mamirubio}J. Martin and M. Milman, \textsl{On the
Calder\'{o}n-Maz'ya-Rubio de Francia extrapolation principle}, preprint 2013.

\bibitem {mamilectures}J. Martin and M. Milman, \textsl{Symmetrization methods
in the theory of Sobolev inequalities}, Lecture Notes, in preparation.

\bibitem {mmp}J. Martin, M. Milman and E. Pustylnik, \textsl{Sobolev
Inequalities: Symmetrization and Self Improvement via truncation}, Journal of
Functional Analysis \textbf{252} (2007), 677-695.

\bibitem {mas}M. Mastylo, \textsl{The Modulus of Smoothness in Metric Spaces
and Related Problems}, Potential Anal. \textbf{35} (2011), 301-328.

\bibitem {maz}V. G. Maz'ya, \textsl{Sobolev Spaces with applications to
elliptic partial differential equations}. Second, revised and augmented
edition. Grundlehren der Mathematischen Wissenschaften [Fundamental Principles
of Mathematical Sciences], \textbf{342}. Springer, Heidelberg, 2011.

\bibitem {maz1}V. G. Maz'ya, \textsl{The} \textsl{p-conductivity and theorems
on imbedding certain functional spaces into a C-space} (Russian), Dokl. Akad.
Nauk SSSR \textbf{140} (1961), 299--302 (English translation: in Soviet Math.
Dokl. 3 (1962).

\bibitem {mie1}E. Milman, \textsl{Concentration and isoperimetry are
equivalent assuming curvature lower bound}, C. R. Math. Acad. Sci. Paris
\textbf{347} (2009), 73--76.

\bibitem {MiE}E. Milman, \textsl{On the role of Convexity in Isoperimetry,
Spectral-Gap and Concentration}, Invent. Math. \textbf{177} (2009), 1-43.

\bibitem {mie2}E. Milman, \textsl{On the role of convexity in functional and
isoperimetric inequalities}, Proc. London Math. Soc., Proc. \textbf{999}
(2009), 32-66.

\bibitem {miemie}E. Milman, \textsl{Isoperimetric and Concentration
Inequalities - Equivalence under Curvature Lower Bound}, Duke Math. J.
\textbf{154} (2010), 207-239.

\bibitem {miemaz}E. Milman, \textsl{A converse to the Maz'ya inequality for
capacities under curvature lower bound}, in A. Laptev (ed), Around the
research of Vladimir Maz'ya I: Function Spaces, Springer, 2010, pp 321-348.

\bibitem {milcon}E. Milman, \textsl{Isoperimetric bounds on convex manifolds},
Contemp. Math \textbf{545} (2011), 195-208.

\bibitem {milocal}M. Milman, \textsl{Local operators vs Lorentz-Marcinkiewicz
spaces}. Interpolation spaces and related topics (Haifa, 1990), Israel Math.
Conf. Proc. 5 (1992), 151--157.

\bibitem {mp}M. Milman and E. Pustylnik,\textsl{ On sharp higher order Sobolev
embeddings}, Comm. Contemp. Math. \textbf{6} (2004), 495-511.

\bibitem {Vmilman}http://en.wikipedia.org/wiki/Vitali\_Milman

\bibitem {okl}E. Oklander, \textsl{Interpolacion, espacios de Lorentz y
teorema de Marcinkiewicz}, Cursos y Seminarios 20, Univ. Buenos Aires, 1965.
(See also E. Oklander, \textsl{On interpolation of Banach spaces}, Thesis,
Univ. Chicago, 1963)

\bibitem {on}R. O'Neil,\textsl{ Convolution operators and L(p,q) spaces}, Duke
Math. J. \textbf{30} (1963), 129--142.

\bibitem {perez}F.J. P\'{e}rez L\'{a}zaro, \textsl{A note on extreme cases of
Sobolev embeddings}, J. Math. Anal. Appl. \textbf{320} (2006), 973--982.

\bibitem {pick}L. Pick, A. Kufner, O. John and S. Fucik, \textsl{Function
Spaces}, Volume 1, Walter de Gruyter \& Co, Berlin, 2012

\bibitem {pis}G. Pisier, \textsl{Factorization of operators through
}$L_{p\infty}$\textsl{ or }$L_{p1}$\textsl{ and non-commutative
generalizations}, Math. Ann. \textbf{276} (1986), 105-136.

\bibitem {pu1}E. Pustylnik, \textsl{On compactness of Sobolev embeddings in
rearrangement-invariant spaces}, Forum Math. \textbf{18} (2006), 839--852.

\bibitem {ra}J. M. Rakotoson, \textsl{R\'{e}arrangement relatif. Un instrument
d'estimations dans les probl\`{e}mes aux limites}, Mathematics \& Applications
\textbf{64}, Springer, Berlin, 2008..

\bibitem {Ros}A. Ros, \textsl{The isoperimetric problem}, In: Global Theory of
Minimal Surfaces. Clay Math. Proc., vol. 2, pp. 175-209, Am. Math. Soc.,
Providence, 2005

\bibitem {rota}G. C. Rota, \textsl{Ten Lessons I wish I had been Taught},
(http://alumni.media.mit.edu/\symbol{126}cahn/life/gian-carlo-rota-10-lessons.html)

\bibitem {saloff coste}L. Saloff-Coste, \textsl{Aspects of Sobolev
inequalities}, Cambridge University Press, 2002.

\bibitem {slavikova}L. Slav\'{\i}kov\'{a}, \textsl{Compactness of higher order
Sobolev embeddings}, Master Thesis, Charles University, Prague 2012.

\bibitem {sut}V. N. Sudakov and B. S. Tsirelson, \textsl{Extremal properties
of half-spaces for spherically invariant measures}. J. Soviet. Math.
\textbf{9} (1978), 918; translated from Zap. Nauch. Sem. L.O.M.I. \textbf{41}
(1974), 1424.

\bibitem {tal}G. Talenti,\textsl{ Inequalities in rearrangement-invariant
function spaces}, in: Nonlinear Analysis, Function Spaces and Applications,
vol. 5, Prometheus, Prague, 1995, pp. 177--230. for a comprehensive bibliography

\bibitem {tar}L. Tartar, \textsl{Imbedding theorems of Sobolev spaces into
Lorentz spaces}, Boll. Unione Mat. Ital. Sez B Artic. Ric. Mat. (8) \textbf{1}
(1998), 479--500.

\bibitem {tr}N. Trudinger, \textsl{On imbeddings into Orlicz spaces and some
applications}, J. Math. Mech. \textbf{17} (1967) 473-483.

\bibitem {tri}H. Triebel, \textsl{Tractable embeddings of Besov spaces into
Zygmund spaces}, Function spaces IX, 361-377, Banach Center Publ. \textbf{92},
Polish Acad. Sci. Inst. Math., Warsaw, 2011.

\bibitem {tri1}H. Triebel, \textsl{Tractable embeddings}, preprint, University
of Jena, Nov. 2012.

\bibitem {xiao}J. Xiao and Z. Zhai, \textsl{Fractional Sobolev,
Moser-Trudinger, Morrey-Sobolev inequalities under Lorentz norms}, J. Math.
Sci. (New York), \textbf{166} (2010), 357-376.
\end{thebibliography}
\end{document}